\documentclass{amsart}
\usepackage{amsthm}
\usepackage{amssymb}
\usepackage{amsmath}

\def	\cal	\mathcal

\def	\bZ	{{\boldsymbol Z}}
\def	\Z	{{\mathbb Z}}
\def	\R	{{\mathbb R}}
\def	\C	{{\mathbb C}}

\def	\t	{{\mathfrak t}}
\def	\g	{{\mathfrak g}}
\def	\ut	{\underline{\mathfrak t}}
\def	\ft	{{\mathfrak t}}
\def	\h	{{\mathfrak h}}

\def	\inv	{^{-1}}

\def	\ss	{\scriptstyle}

\def	\calO	{{\mathcal{O}}}

\def	\conv	{\operatorname{conv}}

\def	\d	{{\operatorname{d}}}

\def	\image	{\operatorname{image}}

\def	\SO	{\operatorname{SO}}
\def	\Gr	{\operatorname{Gr}}
\def	\to	{\longrightarrow}
\def	\ssminus	{\smallsetminus}

\def	\red	{{\operatorname{red}}}

\def	\DH	{Duistermaat-Heckman\ }
\def	\bs	{\boldsymbol}

\def	\half	{\frac{1}{2}}

\def	\tbasic	{{\text{basic}}}
\def	\id	{\text{id}}

\def 	\obas	{{{\Omega^*_\tbasic}}}

\def	\cU		{{\mathcal U}}

\def	\PhiT		{\Phi \mbox{--} T}
\def	\zbar	{{\ol{z}}}

\def	\Fbar	{{\ol{F}}}

\def	\omegabar	{{\tilde{\omega}}}
\def	\i	{{\sqrt{-1} \ }}

\def	\vC	{\text{\v{C}}} 
\def	\vH	{\text{\v{H}}} 

\newcommand     {\ol}[1]   {\overline{#1}}

\newcommand     {\dd}[1]   {\frac{\partial}{\partial #1}}

\newcommand     {\comment}[1]   {}
\newcommand{\mute}[2] {}
\newcommand     {\printname}[1] {}
\newcommand{\labell}[1] {\label{#1}\printname{#1}}

\numberwithin{equation}{section}
\newtheorem {Theorem}			{Theorem}
\newtheorem {refTheorem}[equation]	{Theorem}
\newtheorem {Lemma}[equation]     	{Lemma}  	
  
\newtheorem {Corollary} [equation]	{Corollary}  
\newtheorem {Proposition} [equation]	{Proposition}

\theoremstyle{definition}
\newtheorem{Definition}[equation]{Definition}

\theoremstyle{remark}
\newtheorem{Remark}[equation]{Remark}

\newtheorem{Example}[equation]{Example}

\def	\del	{\partial}

\begin{document}

\newenvironment{pf}{\begin{proof}}{\end{proof}}
\newenvironment{pf*}[1]{\begin{proof}[#1]}{\end{proof}}

\newenvironment{reProposition}[1]
 {\medskip \noindent \textbf{Proposition #1.} 
  \it}
 {\smallskip}

\newenvironment{reTheorem}[1]
 {\medskip \noindent \textbf{Theorem #1.} 
  \it}
 {\smallskip}

\title[Centered complexity one actions]
{Centered complexity one Hamiltonian torus actions}

\author[Yael Karshon]{Yael Karshon}
\address{The Hebrew University of Jerusalem, Giva't Ram, Jerusalem 91904, 
Israel}
\email{karshon@math.huji.ac.il}

\author[Susan Tolman]{Susan Tolman}
\address{Dept.\ of Mathematics, Univ.\ of Illinois, Urbana, IL 61801}
\email{stolman@math.uiuc.edu}

\thanks{Y. Karshon was partially supported by NSF grant DMS-9404404
during earlier work on this project, and by M.S.R.I.\ during the fall
of 1999.
S. Tolman is partially supported by a Sloan fellowship and by 
NSF grant DMS-980305.
Both authors are partially supported by BSF grant 96--00210.}

\begin{abstract}
We consider symplectic manifolds with Hamiltonian torus actions 
which are ``almost but not quite completely integrable": the dimension 
of the torus is one less than half the dimension of the manifold.
We provide a complete set of invariants 
for such spaces when they are ``centered" (see below) 
and the moment map is proper.
In particular, this classifies the moment map preimages of all 
sufficiently small open sets, which is an important step towards
global classification.
As an  application,
we construct a full packing of each of the Grassmanians 
$\Gr^+(2,\R^5)$ and $\Gr^+(2,\R^6)$ by two equal symplectic balls.
\end{abstract}

\maketitle

\tableofcontents

\pagebreak

\section{Introduction}
\labell{sec:intro}

Let a torus $T \cong (S^1)^{\dim T}$ act effectively  on
a symplectic manifold $(M,\omega)$ by symplectic transformations
with a moment map $\Phi \colon M \to \t^*  $, that is,
\begin{equation} \labell{def moment}
        \iota(\xi_M) \omega = - d \left< \Phi,\xi \right>
\end{equation}
for every $\xi$ in the Lie algebra $\t$ of $T$, 
where $\xi_M$ is the corresponding vector field on $M$.
The dimension of the torus is at most half the dimension of the manifold.
The difference $k = \half \dim M - \dim T$ is half the dimension
of the symplectic quotient  $\Phi\inv(\alpha)/T$ 
at a regular value $\alpha \in \Phi(M)$.
We call this number $k$ the \textbf{complexity}\footnote
{
We changed our earlier term \emph{deficiency} to \emph{complexity} 
in order to be consistent with the algebraic geometers' terminology.};
see Definition \ref{def:complexity k}.

The cases when $M$ is compact and the complexity is zero,
also known as \emph{symplectic toric manifolds} or  
\emph{Delzant spaces}, are classified by their moment images \cite{De}.
The first  examples of complexity one spaces are compact symplectic
surfaces (with no action).  By Moser \cite{moser}, 
these are classified by their genus and total area.
Compact symplectic four manifolds with Hamiltonian  circle actions
were classified by the first author \cite{karshon:periodic}; 
also see \cite{ah-ha,audin:paper,audin:book}.
In the algebraic category, complexity one actions (of possibly non-abelian
groups) were recently classified by Timash\"ev \cite{T1,T2}.
Among other works on Lie group actions of complexity zero or one
are \cite{iglesias,Del2,woodward,GSj,knop} in the symplectic category;
\cite{toroidal,or-wag,rynes,FK,b-b,LV} in the algebraic category;
\cite{fintushel,OR} in the smooth category.

This paper is the first in a series of papers in which we 
study complexity one spaces of arbitrary dimension.
In this paper we study the basic building blocks:
the preimages under the moment map of sufficiently small
open subsets in $\t^*$. 
We provide invariants which determine these spaces up to an
equivariant symplectomorphism. 
Moreover, our techniques apply to 
all centered complexity one spaces (see below).

In later papers, we will
give a complete global  classification of complexity one spaces.
This will provide a basis from which to address global questions about
complexity ones spaces, such as
\begin{enumerate}
\item What is the space of automorphisms?
\item When is there a compatible K\"ahler structure?
\item  When are two complexity one spaces symplectomorphic? equivariantly
diffeomorphic?  diffeomorphic?
\end{enumerate}

In this paper, because 
we wish to restrict to the preimages of  open subsets  of $\t^*$,
we do not insist that our manifolds be compact. Instead,
we assume that the moment map is proper as a map to an open convex
set $U \subset \t^*$, that is,
that the preimage of every compact subset of $U$ is compact.
The connectedness and convexity theorems still hold in this generality.
\begin{Definition} \labell{def:complexity k}
Let $T$ be a torus.
A \textbf{proper Hamiltonian $T$-manifold}
is a connected symplectic manifold $(M,\omega)$ 
together with an effective action of $T$,
an open convex subset $U \subseteq \t^*$,
and a proper moment map $ \Phi \colon M \to U$.
Here, $\t$ is the Lie algebra of $T$ and $\t^*$ the dual space.
For brevity, in this paper we call $(M,\omega,\Phi,U)$
a \textbf{complexity $\mathbf k$ space},
where $k = \half \dim M - \dim T$.
An \textbf{isomorphism} between two such spaces over the same set $U$ 
is an equivariant symplectomorphism that respects the moment maps.
\end{Definition}

\begin{Example}
A compact symplectic manifold with a torus action 
and a moment map is a proper Hamiltonian $T$-manifold over $\t^*$.
\end{Example}

\begin{Example}
Let $(M,\omega,\Phi,U)$ be a proper Hamiltonian $T$-manifold.
For any open convex subset $V \subseteq U$,
the preimage $\Phi\inv(V)$ is a proper Hamiltonian $T$-manifold over $V$. 
The fact that it is connected follows 
from the facts that $\Phi \colon \Phi\inv(V) \to V \cap \Phi(M)$ is proper 
and its image and fibers are connected (see Theorem \ref{convexity etc}) 
by easy point-set topology.
\end{Example}

We will now describe invariants of a complexity one space.
Over sufficiently small subsets of $\t^*$, 
these will be enough to characterize the space.

The \textbf{Liouville measure} on a $2n$ dimensional 
symplectic manifold $(M,\omega)$
is given by integration of the volume form $\omega^n / n!$ 
with respect to the symplectic orientation. In the presence of a Hamiltonian
action, the \textbf{\DH measure} is the push-forward of Liouville measure
by the moment map.
It is equal to the \textbf{\DH function} times Lebesgue measure 
on $\t^*$.

Assume $M$ is connected.
For any  value $\alpha \in \Phi(M)$, if the symplectic quotient
$\Phi\inv(\alpha)/T$ is not a single point it 
is homeomorphic to a connected closed oriented surface 
(see Proposition \ref{topology of quotient}). The genus of this
surface does not depend on $\alpha$ (see Corollary \ref{same genus});
we call it the \textbf{genus} of the complexity one space.

The \textbf{stabilizer} of a point $x \in M$ is the closed subgroup
$H = \{\lambda \in T \mid \lambda \cdot x = x\}$. 
The isotropy representation at $x$ is the linear representation
of $H$ on the tangent space $T_xM$. 
Points in the same orbit have the same stabilizer
and their isotropy representations are linearly symplectically isomorphic;
this isomorphism class is the \textbf{isotropy representation} of the
orbit.

An orbit is \textbf{exceptional} if every nearby orbit in the same
moment  fiber has a strictly smaller stabilizer.
Since each moment  fiber is compact,
it contains finitely many exceptional orbits.
The \textbf{isotropy data} at $\alpha \in U$ is the unordered list of
isotropy representations of the exceptional orbits in $\Phi\inv(\alpha)$.

With these definitions on hand, let us state our 
main theorem, which gives necessary and sufficient conditions
for two complexity one spaces to be locally isomorphic.

\begin{Theorem}[Local Uniqueness] \labell{local uniqueness}
Let $(M,\omega,\Phi,U)$ and $(M',\omega',\Phi',U)$ be 
complexity one spaces.
Assume that their \DH measures are the same and that their genus 
and isotropy data over $\alpha \in U$ are the same. 
Then there exists a neighborhood of $\alpha$ 
over which the spaces are isomorphic.
\end{Theorem}

Here is a simple proof of Theorem \ref{local uniqueness} 
in the case that the torus action on $\Phi\inv(\alpha)$ is free:

\begin{quotation}
The symplectic quotient $\Phi\inv(\alpha) / T$ is a symplectic surface.
Its symplectic area is the value of the \DH function at $\alpha$.
Together with the genus, this determine the surface.

The moment fiber $Z : = \Phi\inv(\alpha)$ is a principal $T$-bundle 
over the symplectic quotient.
Its Chern class is given by the slope of the \DH function at $\alpha$ \cite{DH}.

The pullback to the moment fiber $Z$ of the symplectic form on the symplectic
quotient is the restriction $i_Z^* \omega$ of $\omega$ to $Z$.
By the equivariant coisotropic embedding theorem 
(see \cite[lecture 5]{w:lectures}).
a neighborhood of the moment fiber $Z$ is determined up to 
equivariant symplectomorphism by $(Z,i_Z^*\omega)$.
Since the moment map is proper,
this neighborhood contains the preimage of a neighborhood of $\alpha$.
\end{quotation}

This argument straightforwardly
extends to the case that $\alpha$ is a regular 
value of the moment map.
The main volume of this paper consists of carefully
extending the argument to singular values of the moment map.

Additionally, we prove a variation of the theorem 
for ``centered spaces" over larger subsets of $\t^*$. 
Recall that the \textbf{orbit type strata} are the connected 
components of the sets of points with the same stabilizer.

\begin{Definition} \labell{centered-definition}
A proper Hamiltonian $T$-manifold $(M,\omega,\Phi,U)$ 
is \textbf{centered} about a point $\alpha \in U$ if
$\alpha$ is contained in the closure of the moment image
of every orbit type stratum in $M$.
\end{Definition}

\begin{Theorem}[Centered Uniqueness] \labell{thm:centered-uniqueness}
Let $(M,\omega,\Phi,U)$ and $(M',\omega',\Phi',U)$ be 
complexity one spaces that are centered about $\alpha \in U$.
Assume that their \DH measures are the same and that their genus and 
isotropy data over $\alpha \in \t^*$ are the same. Then the spaces are 
isomorphic. 
\end{Theorem}

Finally, we present an application of our results to symplectic topology:
we construct full packings of two Grassmanians.
Holomorphic techniques in symplectic topology are useful in 
giving obstructions to embeddings, but have been less successful in 
constructing embeddings.   Our construction, like several
previous ones \cite{traynor,karshon:appendix}.
uses equivariant techniques to solve this non-equivariant problem.
In a future paper \cite{gromov}, we will extend the techniques that
we develop here to address this question more deeply. 
Here, we content ourselves with a simple, but fairly representative, 
application:

\begin{Theorem} \labell{thm:applications}
Let $M$ be the Grassmanian $\Gr^+(2, \R^5)$ or $\Gr^+(2, \R^6)$. 
There exists an equivariant symplectic embedding of a disjoint union of
two symplectic balls with linear actions and with equal radii into $M$
such that the complement of the image has zero volume.
A fortiori, each of these Grassmanians 
can be fully packed by two equal symplectic balls.
\end{Theorem}

\subsubsection*{Acknowledgement}
We thank F. Knop and D. Luna for explaining to us 
aspects of complexity one actions related to algebraic geometry.

\section{Background}
\labell{sec:background}

We now set our notation and review some background material.

Let a torus $T$ act effectively on a symplectic manifold $(M,\omega)$.
The symplectic slice at $x \in M$ is the symplectic vector space
$$(T_x\calO)^\omega/(T_x\calO \cap (T_x\calO)^\omega),$$ 
where $\calO$ is the $T$-orbit of $x$ in $M$.
Let $H \subset T$ be the stabilizer of $x$.
The isotropy representation of $H$ on $T_xM$ induces a representation
on the symplectic slice, called the \textbf{slice representation}.
The slice representation is isomorphic to the action of $H$ on $\C^n$ 
through an inclusion  $\rho = (\rho_1 , \ldots, \rho_n) \colon H \to (S^1)^n$.
The isotropy representation is the direct sum
of the slice representation and a trivial representation.
The isotropy characters $\rho_i$ are determined up to permutation.
The differential of each $\rho_i \colon H \to S^1$ is an element $\eta_i$
of the dual space $\h^*$.
The $\eta_i$ are called the \textbf{isotropy weights}.

We fix an inner product on the Lie algebra $\t$ of our torus $T$, 
once and for all.
This determines a projection $\t \to \h$ and, dually, an inclusion
$\h^* \hookrightarrow \t^*$ for any subspace $\h \subset \t$.
Throughout this paper, we will identify $\h^*$ with its image in $\t^*$.

The Guillemin-Sternberg-Marle local normal form theorem 
classifies the neighborhoods of orbits in symplectic manifolds 
with Hamiltonian actions of compact groups. We state it for tori:

\begin{refTheorem}[Local normal form] 
Let a closed subgroup $H$ of a torus $T$ act on $\C^n$ 
by an inclusion $\rho \colon H \to (S^1)^n$ with weights
$\eta_1,\ldots,\eta_n$ and moment map
\begin{equation*} 
\Phi_H(z) = \frac{1}{2} \sum_{j=1}^n |z_j|^2 \eta_j. 
\end{equation*}
\begin{enumerate}
\item
Equip $T^*(T) \times \C^n$ 
with the standard symplectic form and the diagonal $H$ action.
Its symplectic quotient by $H$  
can be  identified with the model
\begin{equation*} 
Y = T \times_H \C^n \times \h^0.
\end{equation*}
Given $\alpha \in \t^*$,
\begin{equation*} 
 \Phi_Y ( [t,z,\nu]) =  \alpha + \Phi_H(z) + \nu 
\end{equation*}
is  a moment map for the left $T$ action.
Here,  $T^*(T) = T \times \t^*$ is the cotangent bundle of $T$
\item
Let the torus $T$ act effectively on a symplectic manifold  $(M,\omega)$  
with a moment map $\Phi \colon M \to \t^*$.
Given a point  $x \in M$ with slice representation $\rho$,
there exists a neighborhood of the orbit $T \cdot x$ that is
equivariantly symplectomorphic to a neighborhood of the
orbit $\{[t,0,0]\}$  in the model $Y$ with $\alpha = \Phi(x)$.
\end{enumerate}
\end{refTheorem}

This is proved in \cite{GS:normal} and in \cite{marle}.

The local normal form theorem implies an important special case
of Theorem \ref{local uniqueness}:

\begin{Proposition} \labell{proper}
Let $(M,\omega,\Phi,U)$ and $(M',\omega',\Phi',U)$ 
be proper Hamiltonian $T$-manifolds.
Consider a value $\alpha \in U$ so that the moment fibers
$\Phi\inv(\alpha)$ and $\Phi'{}\inv(\alpha)$ each consists of a single orbit.
Suppose that these orbits have the same slice representation. 
Then there exists a neighborhood $V$ of $\alpha$ over which $M$ and $M'$ 
are isomorphic.
\end{Proposition}

\begin{proof}
Denote the orbits over $\alpha$ in $M$ and $M'$ by $\calO$ and $\calO'$.
A neighborhood of $\calO$ in $M$ and a neighborhood of $\calO'$ in $M'$
are each isomorphic to 
a neighborhood of $\{[t,0,0]\}$ in the same local model $Y$, 
by assumption and by the local normal form theorem.
Hence a neighborhood $W$ of $\calO$ is isomorphic to a neighborhood $W'$
of $\calO'$.  Since the moment maps are proper, 
if $V \subseteq U$ is a small enough neighborhood of $\alpha$, 
the preimages $\Phi\inv(V)$ and $\Phi'{}\inv(V)$
are contained in $W$ and $W'$, and are hence isomorphic.
\end{proof}

We will also use the following global properties:

\begin{refTheorem} \labell{convexity etc}
Every proper Hamiltonian $T$-manifold $(M,\omega,\Phi,U)$ 
has the following properties.
\begin{description}
\item[Convexity]
The moment image, $\Phi(M)$, is convex.
\item[Connectedness]
The moment fiber, $\Phi\inv(\alpha)$, is connected for all $\alpha \in U$.
\item[Stability]
As a map to $\Phi(M)$, the moment map is open. 
\end{description}
\end{refTheorem}

These properties, which are intimately related to each other,
are due to Atiyah, Guillemin and Sternberg in the compact case.
For Convexity and Connectedness, see \cite{atiyah}.
For Convexity and Stability, see \cite{GS:convexity}.  
For proper moment maps to open convex sets and a brief history, 
see \cite{LMTW}.

\section{Eliminating the symplectic form}
\labell{sec:eliminate-omega}

Our next task is to  free ourselves from the symplectic form. 
In this section we show that, instead
of working with equivariant symplectomorphisms, it is enough to work
with equivariant diffeomorphisms that respect the orientation and the
moment map.
These are much easier to work with, as one can  apply
techniques from differential topology.

\begin{Definition}
Let $M$ and $M'$ be oriented manifolds with  $T$ actions and
$T$-invariant maps $\Phi \colon M \to \t^*$ and $\Phi' \colon M' \to \t^*$.
A {\bf $\bs \PhiT$--diffeomorphism} 
from $(M,\Phi)$ to  $(M,\Phi')$ is an orientation preserving
equivariant diffeomorphism $\Psi \colon M \to M'$ that satisfies $\Psi^*(\Phi')
 = \Phi$.
\end{Definition}

In this section and the next one we will need to assume the 
following technical condition:
\begin{equation} \labell{technical121}
\begin{array}{l}
\text{The restriction map }
	H^2(M/T,\Z) \to H^2(\Phi\inv(y)/T,\Z) \\
\text{is one-to-one for some regular value $y$ of $\Phi$.}
\end{array}
\end{equation}
In fact, this restriction map is always one-to-one.  
We prove this for the preimage of small open sets 
in this paper; see Corollary \ref{cor 121} 
and Lemma \ref{121 near point}.
In a later paper, we will prove it for all complexity one spaces. 

\begin{Proposition} \labell{eliminate}
Let $(M,\omega,\Phi,U)$ and $(M',\omega',\Phi',U)$ be 
complexity one spaces that satisfy Condition \eqref{technical121}.
Assume that 
they have the same  \DH measure.\footnote{ \labell{DH constant}
In fact, we only need that the \DH functions agree at a point. Contrast
with footnote \ref{DH slope}.}
Then there exists an equivariant symplectomorphism from $M$ to $M'$  
if and only if 
there exists a $\PhiT$--diffeomorphism from $M$ to $M'$.
\end{Proposition}

The following proof of Proposition \ref{eliminate} relies on a couple
of technical lemmas which we postpone until after the proof.

\begin{proof}
Let $g \colon M \to M'$ be  a $\PhiT$--diffeomorphism.
By  Lemma \ref{cohomology}, there exists 
a basic  one-form $\beta$  on $M$ (see Remark~\ref{basic})
such that $d\beta = g^*\omega' - \omega$.
We now apply Moser's method:

Define $\omega_t := (1-t) \omega + t g^* \omega'$ 
for all $0 \leq t \leq 1$. 
By Lemma \ref{nondegenerate} below,
the $\omega_t$ are nondegenerate. 
Let $X_t$ be the vector field determined by $i_{X_t} \omega_t = -\beta$.
The vector field $X_t$ preserves the level sets of $\Phi$,
because for every $\xi \in \t$,
$\left< d\Phi(X_t), \xi \right>  
= - \omega_t(\xi_M , X_t) = - i_{\xi_M} \beta = 0$.
Since $\Phi$ is proper, the time-dependent vector-field $X_t$ 
integrates to a flow, $F_t \colon M \to M$.
Let $g_t = g \circ F_t$.  Then $\Phi'  \circ g_t = \Phi$ 

The vector field $X_t$ is invariant 
because $\omega_t$ and $\beta$ are invariant. 
Consequently, $g_t$ is  equivariant,
and hence is a $\PhiT$--diffeomorphism.  Finally,
\begin{eqnarray*}
	\frac{d}{dt}(F_t^*\omega_t) & = &
	F^*_t( L_{X_t} \omega_t) + F^*_t (\frac{d}{dt} \omega_t) \\
 	& = & F^*_t d \iota_{X_t} \omega_t + F^*_t (\omega_1 - \omega_0) \\
	& = & F^*_t( -d\beta + d\beta) = 0.
\end{eqnarray*}
Therefore, $F_1^* \omega_1 = F_0^* \omega_0 = \omega_0$,
since $F_0$ is the identity.
Then,  $g_1^* \omega' = F_1^* (g^* \omega') =  
F_1^*(\omega_1) = \omega_0 = \omega$.
In other words,  $g_1$ is an (equivariant) symplectomorphism.
\end{proof}

Before proving the technical lemmas used in the above proof,
let us recall the notion of basic forms:

\begin{Remark} \labell{basic}
Let a compact Lie group $G$ act on a manifold $M$,
and for $\xi \in \g$ let $\xi_M$ be the generating vector-fields.
A differential form $\beta$ on $M$ is \textbf{basic} if
it is $G$ invariant and horizontal, that is,
$\iota_{\xi_M} \beta =0 $ for all $\xi \in \g$.
The basic differential forms on $M$ constitute a differential
complex $\obas(M)$ whose cohomology coincides with the \v{C}ech cohomology
of the topological quotient, $M/G$.  See \cite{koszul}. 
To see this, repeat the standard \v{C}ech-de Rham
spectral-sequence argument, as in \cite{BT}. 
It still works because,
by  the local normal form for smooth actions of
compact Lie groups, every orbit in $M$
has a neighborhood on which the complex of basic forms is acyclic.
\end{Remark}

\begin{Lemma} \labell{cohomology} 
Let $(M,\omega,\Phi,U)$ and $(M',\omega',\Phi',U)$ be 
complexity one spaces.
Assume that  the \DH measures for $M$ and  $M'$ are the same, and
that the spaces satisfy Condition \eqref{technical121}.
Then for every $\PhiT$--diffeomorphism $g \colon M \to M'$  
there exists  a basic one-form $\beta$ on $M$
such that $d \beta = g^* \omega' - \omega$.
\end{Lemma}

\begin{proof}
Let $\Omega = g^* \omega' - \omega$.
Since $\omega$ and $g^* \omega'$ are closed invariant symplectic forms 
on $M$ with the same moment map,
$\iota(\xi_M) \Omega =0$ for all $\xi \in \t$.
Since $\Omega$ is  also invariant, it is basic.

By Condition \eqref{technical121} it suffices to show 
that the restriction of $\Omega$ to the fiber $\Phi\inv(\alpha)$ 
is exact for some regular value $\alpha$ of $\Phi$.
Since this restriction is the pull-back of a differential form 
$\Omega_\red$ on the orbifold $M_\red = \Phi\inv(\alpha) / T$, 
it is enough to show that $\Omega_\red$ is exact.
Since $M_\red$ is two dimensional,
it is enough to show that the integral of $\Omega_\red$ over it is zero,
i.e., that the integrals of $\omega$ and $g^*\omega'$ are equal.
But this follows from the fact that the \DH measures for $M$ and $M'$
are the same, because the density functions for these measures are
given by the symplectic volumes of the symplectic quotients; 
see \cite[\S 3]{DH}.
\end{proof}

\begin{Lemma} \labell{nondegenerate}
Let $(M,\omega,\Phi,U)$ and $(M',\omega',\Phi',U)$ be complexity one spaces.
Let $g \colon M \to M'$ be  a $\PhiT$--diffeomorphism. 
Then the two-form $\omega_t = (1-t)\omega + t g^*\omega'$ 
is nondegenerate for all $0 \leq t \leq 1$.
\end{Lemma}

\begin{proof}
First, let a compact abelian Lie group $H$ act on $\C^n$ as a codimension one
subgroup of $(S^1)^n$ with isotropy weights 
$\eta_1,\ldots,\eta_n \in \h^*$.
The vector fields for this action are
$$
	\xi_M = \sqrt{-1}\ \sum_i \eta_i (\xi)
	\left( z_i  \dd{z_i} - \zbar_i  \dd{\zbar_i} \right), 
	\quad \xi \in \h.
$$

Let $\omegabar_0$ and $\omegabar_1$ be invariant symplectic forms on $\C^n$
with constant coefficients that have the same moment map, $\Phi_H$, 
and that induce the same orientation.  We will show that 
$\omegabar_t = (1-t) \omegabar_0 + t \omegabar_1$
is non-degenerate.

Because $\omegabar_t$ is real valued, it can be written in the form
$$ \begin{array}{lll}
	\omegabar_t & = & \sqrt{-1} \sum\limits_j A_j^t dz_j \wedge d\overline{z_j} \\
	& + & \half \sqrt{-1} \sum\limits_{j \neq k} \left( 
	      B_{jk}^t dz_j \wedge dz_k  
		  - \overline{ B_{jk}^t}  d\overline{z_j} \wedge d\overline{z_k} 
		  \right)  \\
	& + & \sqrt{-1} \sum\limits_{j\neq k} C_{jk}^t dz_j \wedge d\overline{z_k},
\end{array}
$$
where $A_j^t$ are real, $B_{jk}^t$ and $C_{jk}^t$ are complex,
$B_{jk}^t = -B_{kj}^t$, and $C_{jk}^t = \ol{C_{kj}^t}$.

By the definition of moment map,
$\del \Phi_H^\xi / \del z_j$ is the coefficient 
of $dz_j$ in $-\iota(\xi_M) \tilde{\omega}_t$. Hence 
$$ \del \Phi_H^\xi / \del z_j = 
  \eta_j(\xi) A_j^t \zbar_j
  + \sum\limits_{k \neq j} \eta_k(\xi) 
    \left( B_{kj}^t z_k + C_{jk}^t \zbar_k \right). $$
Differentiating again, 
$\del^2 \Phi_H^\xi / \del z_j \del \zbar_j = \eta_j(\xi) A_j^t$.
Because $\omegabar_t$ is $T$-invariant, $B_{jk}^t = 0$ unless 
$\eta_j = - \eta_k$. In this case,
$\del^2 \Phi_H^\xi / \del z_j \del z_k = \eta_k(\xi) B_{kj}^t.$
Finally,   $C_{jk}^t = 0$ unless $\eta_j = \eta_k$;
in this case,
$ \del^2 \Phi_H^\xi / \del z_j \del \overline{z_k} 
  = \eta_k(\xi) C_{jk}^t$.
Therefore, $\eta_j(\xi) A_j^t$, $\eta_k(\xi) B_{kj}^t$, 
and $\eta_k(\xi) C_{jk}^t$ are determined by $\Phi_H^\xi$.

By what we have shown, if $\eta_j \neq 0$, 
the coefficients $A_j^t$, $B_{jk}^t$, and $C_{jk}^t$ 
are determined by $\Phi_H$.  Thus, if no weight is zero, 
$\omegabar_t$ is independent of $t$, hence it is non-degenerate.
Since the action is effective and the dimension of $\h^*$ is $n-1$, 
the only other possibility is that exactly one of the weights -- 
let's say the first one -- is zero,
and the others form a basis of $\h^*$.
In this case, $B_{ij}^t = C_{ij}^t = 0$ for all $i$ and $j$,
and so the top power of $\tilde{\omega}_t$ 
is $\prod_{j=1}^n A_j^t$ times the standard volume form. 
Since $\omegabar_0$ and $\omegabar_1$ induce the same orientation,
and since $A_j^t$ is determined by $\Phi_H$ and hence independent of $t$
for $j \neq 1$, the signs of $A_1^0$ and $A_1^1$ are the same. 
Therefore, $A_j^t = (1-t) A_j^0 + t A_j^1$ is never zero, 
and $\omegabar_t$ is non-degenerate.

Now let $x \in M$ be any point with stabilizer $H \subset T$.
By the local normal form theorem,
a neighborhood of the orbit $T \cdot x$ in $M$ 
with the symplectic form $\omega_0$
is equivariantly symplectomorphic
to a neighborhood of the orbit $\{ [t,0,0] \}$ in the model
$T \times_H \C^n \times \h^0$.
The tangent space at $x$ splits as
$$T_x M = \t / \h \oplus \h^0 \oplus \C^n$$ 
where $\t/\h$ is the tangent to the orbit.
By the definition of the moment map, $\omega_t|_x$ is given by a
block matrix of the form
$$ \left( \begin{array}{ccc} 0 & I & 0 \\ -I & * & * \\ 
	0 & * & \omegabar_t \end{array} \right)
$$
where $I$ is the natural pairing between the vector space $\t/\h$
and its dual, $\h^0$, and where $\omegabar_0$ and $\omegabar_1$ are
linear symplectic forms on $\C^n$ with the same moment map and the 
same orientation.
By the above argument, $\omegabar_t$ is nondegenerate. 
Consequently, $\omega_t|_x$ is nondegenerate.
\end{proof}

\section{Passing to the quotient}
\labell{sec:quotient}

In this section we show that, as long as two complexity one spaces have 
the same \DH measure, we can  reduce the problem of finding
a $\PhiT$--diffeomorphism between them to the easier problem
of finding a $\Phi$-diffeomorphism between their quotients. 
Some techniques in this section are adapted from Haefliger and Salem
\cite{HS}.

Let a compact torus $T$ act on a manifold $N$.
The quotient $N/T$ can be given
the quotient topology and a natural differential structure,
consisting of the sheaf of real-valued functions
whose pullbacks to $N$ are smooth. 
We say that a map $h \colon N/T \to N'/T$ is smooth if it pulls back smooth 
functions to smooth functions; it is a diffeomorphism if it is smooth 
and has a smooth inverse. See \cite{schwarz:IHES}.
If $N$ and $N'$ are oriented, the choice of an orientation on $T$ 
determines orientations on the smooth part of $N/T$ and $N'/T$.  
Whether or not a diffeomorphism $f\colon N/T \to N'/T$ preserves orientation 
is independent of this choice. 

While this notion of diffeomorphism is natural, 
we will also need another less natural but stronger notion.

\begin{Definition} \labell{Phi diffeo}
Let $M$ and $M'$ be oriented manifolds with  $T$ actions and
$T$-invariant maps $\Phi\colon M \to \t^*$ and $\Phi' \colon M' \to \t^*$.
A $\mathbf{\Phi}$\textbf{-diffeomorphism}  from $M/T$ to $M'/T$
is an orientation preserving diffeomorphism  $\Psi \colon M/T \to M'/T$
such that
\begin{enumerate}
\item $\Psi$ preserves the moment map, i.e., $\Psi^* \Phi' = \Phi$.
\item
Each of $\Psi$ and $\Psi\inv$ lifts to a $\PhiT$--diffeomorphism
in a neighborhood of each exceptional orbit. 
\end{enumerate}
\end{Definition}

\begin{Proposition} \labell{prop:HS}
Let $(M,\omega,\Phi,U)$ and $(M',\omega',\Phi',U)$ be 
complexity one spaces.
Assume that Condition \eqref{technical121} is satisfied.
If both spaces have the same \DH measure,\footnote{ 
\labell{DH slope} 
In fact, we only need their \DH functions to have the same slope; 
they may differ by a constant.  Contrast with footnote \ref{DH constant}.} 
every $\Phi$-diffeomorphism from $M/T$ to $M'/T$
lifts to $\PhiT$--diffeomorphism from $M$ to $M'$.
\end{Proposition}

The first step in proving this proposition is to show that
on the non-exceptional orbits every $\Phi$-diffeomorphism lifts locally 
to a $\PhiT$--diffeomorphism. We do this in the next three lemmas.

\begin{Lemma}
\labell{non-exceptional}
Every local model for a non-exceptional orbit 
in a complexity one space has the form
\begin{equation} \labell{non-exceptional model}
 Y = T \times_H \C^h \times \C \times \h^0,
\end{equation}
where $H \subseteq T$ is a closed $h$ dimensional subgroup 
which acts on $\C^h$ through an isomorphism with $(S^1)^h$.
\end{Lemma}

\begin{proof}
Let $Y := T \times_H \C^{h+1} \times \h^0$ be the local model
for a non-exceptional orbit in $\Phi\inv(\alpha)$, with $h=\dim H$.
Inside the moment fiber $\Phi_Y\inv(\alpha)$,  
the set of points with stabilizer $H$ is 
\begin{equation} \labell{stratum}
T \times_H (\C^{h+1})^H \times \{ 0 \} ,
\end{equation}
where $(\C^{h+1})^H$ is the subspace fixed by $H$.    
By the definition of exceptional orbit, this subspace is not trivial.
Therefore, the local model becomes \eqref{non-exceptional model},
where the group $H$ acts trivially on $\C$ and acts on $\C^h$ 
through an inclusion into $(S^1)^h$.
By a dimension count, this inclusion must be an isomorphism.
\end{proof}

The following lemma tells us that neighborhoods of nonexceptional orbits
can be read off from the moment image.

\begin{Lemma} \labell{nonexceptional near boundary}
Let $(M,\omega,\Phi,U)$ be a complexity one space.
Assume that the  preimage  $\Phi\inv(\alpha)$ of $\alpha \in U$ 
contains a non-exceptional orbit.

There exists a closed connected subgroup $H \subseteq T$ with Lie algebra $\h$
and a basis $\{ \eta_j\}$ for the weight lattice in $\h^*$ so that 
\begin{enumerate}
\item
The group $H$ is the stabilizer 
and the $\eta_j$ are the non-zero isotropy weights of every
non-exceptional orbit in $\Phi\inv(\alpha)$.
\item
In a neighborhood of $\alpha$, the image $\Phi(M)$ 
coincides with the Delzant cone $\alpha + \h^0 + \sum_j \R_+ \eta_j$. 
\end{enumerate}
\end{Lemma}

\begin{proof}
Consider the slice representation at any non-exceptional orbit.
By Lemma \ref{non-exceptional} above, the stabilizer $H$ is connected, 
Since the action is effective, the isotropy weights $\eta_j$
generate the weight lattice.
The stabilizer and these weights 
are determined by the image of the moment map for a local model;
this image is the Delzant cone $\h^0 + \sum_j \R_+ \eta_j$.
Finally, by the stability of the moment map, the image of the moment map
is the same for every local model in $\Phi\inv(\alpha)$.
\end{proof}

\begin{Corollary} \labell{nonex=free}
Over the interior of the moment image, the nonexceptional orbits
are precisely the free orbits.
\end{Corollary}

\begin{Lemma} \labell{local lift on Y}
Let $Y$ be a local model for a non-exceptional orbit with a moment map
$\Phi_Y \colon Y \to \t^*$.  
Let $W$ and $W'$ be invariant open subsets of $Y$.
Let $g \colon W/T \to W'/T$  be a diffeomorphism which preserves 
the moment map. 
Then $g$ lifts to an equivariant diffeomorphism from $W$ to $W'$.
\end{Lemma}

\begin{proof}
Assume $W=W'=Y$; the general case is similar.
Since, by Lemma \ref{non-exceptional}, 
$Y = T \times_H \C^h \times \C \times \h^0$,
we can identify $Y/T$ with $\h^0 \times (\C^h/H) \times \C$.
Since $g$ preserves the moment map, it necessarily has the form
$$ g (\nu, [z] , \zeta) = (\nu, [z], \psi(\nu,[z],\zeta)),$$
for some $\psi \colon \h^0 \times (\C^h/H) \times \C \to \C$.
Similarly, its inverse sends $(\nu, [z], \zeta)$ to
$(\nu, [z], \gamma(\nu,[z],\zeta))$, 
where $\zeta \colon \h^o \times (\C^h/H) \times \C \to \C$.
Since both $g$ and its inverse are smooth,  $\zeta$ and $\gamma$ must
themselves be smooth.  

We define $\tilde{g} \colon Y \to Y$ by
$\tilde{g}([t, z, \zeta, \nu])  = [t,z,\psi(\nu,[z],\zeta), \nu]$.
Then $\tilde{g}$ is a smooth equivariant lift of $g$,
and it has a smooth inverse given by
$[t, z, \zeta, \nu] \mapsto [t,z,\gamma(\nu,[z],\zeta), \nu]$.
\end{proof}

We deduce that a $\Phi$-diffeomorphism lifts to a $\PhiT$--diffeomorphism
locally; we still need to show that, in the proper circumstances,
a $\Phi$-diffeomorphism that lifts locally also lifts globally.
We do this in the lemma below;
the basic idea is that the \DH measure determines the ``fibration" $M \to M/T$.

\begin{Lemma} \labell{Lemma:HS}
Let $(M,\omega,\Phi,U)$ and $(M',\omega',\Phi',U)$ be complexity one spaces.
Assume that Condition \eqref{technical121} is satisfied.
If both spaces have the same \DH measure, 
every homeomorphism from $M/T$ to $M'/T$ that locally lifts
to a $\PhiT$--diffeomorphism also lifts globally 
to a $\PhiT$--diffeomorphism from $M$ to $M'$.
\end{Lemma}

\begin{pf}
Let $T,\t$, and $\ell$ denote the sheaves of 
smooth functions from $M/T$ to $T$, $\t$, and $\ell$, respectively.
Here, $\ell$ denotes the lattice in $\t$.
Let $\ut$ denote the sheaf of locally constant function to $\t$. 

Fix a homeomorphism $\Psi \colon M/T \to M'/T$ that lifts locally.
Choose a cover $\cU$ of $M$ by open invariant sets,
and on each $U_i \in \cU$ a $\PhiT$--diffeomorphism $\Psi_i \colon U_i \to M'$
that is a lift of $\Psi$.
By Lemma \ref{schwartz}  below,
there exist smooth invariant functions $g_{ij} \colon U_i \cap U_j \to T$
such that $g_{ij} \cdot \Psi_j = \Psi_i$ for all $i$ and $j$.
These functions form a \v{C}ech cocycle $g \in \vC^1(\cU,T)$.
The map $\Psi$ will lift to a global $\PhiT$--diffeomorphism exactly 
if the corresponding cohomology class $[g] \in \vH^1(M/T,T)$ is trivial.

The short exact sequence $0 \to \ell \to \ft \to T \to 0$
induces a long exact sequence in cohomology.
Since there exists a smooth partition of unity on $M/T$, 
the cohomology $\vH^i(M/T,\t)$ vanishes for all $i > 0$.
Therefore, $\vH^1(M/T,T) =  \vH^2(M/T,\ell)$.
Condition \eqref{technical121} implies that
the restriction map $\vH^2(M/T,\ell) \to \vH^2(\Sigma,\ell)$ is one-to-one,
where $\Sigma = \Phi\inv(\alpha)/T$ is a regular symplectic quotient. 
Therefore, it is enough to show that the image of
$[g]$ in $\vH^2(\Sigma,\ell)$ is zero.
Since $\vH^2(\Sigma,\ell)$ is torsion free,
it is enough to show that the image of $[g]$ in $\vH^2(\Sigma,\ut)$ vanishes. 
The \v{C}ech-de Rham isomorphism for basic forms on $\Phi\inv(\alpha)$
(see Remark \ref{basic})
takes this image to the cohomology class of the basic differential two-form 
whose restriction to each open set $U_i \cap \Phi\inv(\alpha)$ is
\begin{equation} \labell{CdR}
\pm \sum\limits_{j} d \lambda_j g_{ij}^{-1} dg_{ij},
\end{equation}
(the sign depending on conventions),
where $\{\lambda_i\}$ is a partition of unity subordinate to 
$\cU \cap \Phi\inv(\alpha)$.
We claim that this is exact as a basic form.

Let $\Theta'$ be a connection one-form on ${\Phi'}\inv(\alpha) \subset M'$,
that is, a $T$-invariant $\t$-valued one-form
such that $\Theta'(\xi_{M'}) \equiv \xi$ for all $\xi \in \t$.
Then $\Theta = \sum \lambda_i {\Psi_i}^* \Theta'$ is a connection
one-form on $\Phi\inv(\alpha) \subset M$.
The curvature forms $d\Theta$ and $d\Theta'$ are basic.
Their integrals over the symplectic quotients
are equal to the slopes of the \DH function of $M$ and of $M'$ 
at $\alpha$ \cite{DH}. 
Since these slopes are the same, and since $\Phi\inv(\alpha)/T$
is a two dimensional orbifold, the difference between $d\Theta$
and $\Psi^* d\Theta'$ is exact as a basic form. 
A simple computation shows that this difference is equal to \eqref{CdR}.
\end{pf}

In the above proof we used the following theorem 
of Haefliger and Salem, based on a lemma of Schwarz.

\begin{refTheorem}[\cite{HS}]\labell{schwartz}
Let a torus $T$ act on a manifold $M$. Let $h \colon M \to M$ be an
equivariant diffeomorphism that sends each orbit to itself.
Then there exists a smooth invariant function $f\colon M \to T$
such that $h(m) = f(m) \cdot m$ for all $m \in M$.
\end{refTheorem}

We are finally ready to prove our main proposition.

\begin{proof}[Proof of Proposition \ref{prop:HS}]
Let $\calO$ be a non-exceptional orbit in $M$.
By Definition \ref{Phi diffeo}, any $\Phi$-diffeomorphism sends it
to a non-exceptional orbit $\calO'$ in $M'$.
By Lemma \ref{nonexceptional near boundary}, the local models
for $\calO$ and $\calO'$ are the same $Y$. By Lemma \ref{local lift on Y}
and the local normal form theorem, the map lifts to a $\PhiT$--diffeomorphism
from a neighborhood of $\calO$ to a neighborhood of $\calO'$.
By Lemma \ref{Lemma:HS}, it lifts globally.
\end{proof}

\section{Symplectic representations}
\labell{sec:repr}

So far, we have shown  that two complexity one spaces with the same
\DH measure are isomorphic if their quotients are $\Phi$-diffeomorphic
as long as the spaces satisfy Condition \eqref{technical121}.
The rest of the paper is dedicated to proving that, over small subsets
of $\t^*$, the quotients are indeed $\Phi$-diffeomorphic if the spaces have the 
same genus and isotropy data, and that Condition \eqref{technical121} is always 
satisfied.

In preparation for this, in this section we analyze
how the weights of a symplectic representation
can be read from  its moment map.
The key ingredient, which we use repeatedly,
is simply the formula for the moment map:
let a compact abelian group $H$ act effectively on 
$\C^n$ as a subgroup of $(S^1)^n$ with
weights $\eta_1,\ldots,\eta_n$.  Then
\begin{equation} \labell{PhiH again}
 \Phi_H(z) = \half \sum_{j=1}^n |z_j|^2 \eta_j 
 \end{equation}
is a moment map.

\begin{Lemma}  \labell{positive}
Let a compact abelian group $H$ act effectively on $\C^n$ 
with weights $\eta_1, \ldots, \eta_n$.
The moment map $\Phi_H$ is onto if and only if there exist 
$\xi_j > 0$ so that $\sum \xi_j \eta_j = 0$.
\end{Lemma}

\begin{proof}
Suppose that the moment map $\Phi_H$ is onto. Then every element of $\h^*$
is in the non-negative span of the $\{ \eta_j \}$. 
In particular,  there  exist $a_j \geq 0$ such that 
$ \sum a_j \eta_j = \sum_j -\eta_j$, that is, 
$\sum (1+a_j) \eta_j = 0$.  Let $\xi_j = 1 + a_j$. 

Conversely, suppose that there exist positive
$\xi_j$'s so that $\sum \xi_j \eta_j = 0$. 
Let $\alpha \in \h^*$ be any element.
Because the action is effective, its weights, $\eta_j$, span $\h^*$,
so there exist $a_j$, $j=1,\ldots,n$, such that
$\alpha = \sum_j a_j \eta_j$. Because $\sum_j \xi_j \eta_j =0$,
we also have
$\alpha = \sum_j (a_j + t\xi_j) \eta_j$ for any $t \in \R$.
Because $\xi_j>0$ for all $j$, if we take $t$ large enough we get
that $\alpha$ is in the positive span of the $\eta_j$.  
\end{proof}

\begin{Lemma} \labell{criterion not proper}
Let a compact abelian group $H$ act effectively on $\C^n$ 
with moment map $\Phi_H$.
The moment map $\Phi_H$ is not proper if and only if 
there exist $\xi_j \geq 0$, not all zero,
such that $\sum \xi_j \eta_j = 0$.
\end{Lemma}

\begin{proof}
Suppose that $\sum \xi_j \eta_j = 0$ for some $\xi_j \geq 0$,
not all zero.
Since  the moment fiber $\Phi_H\inv(0)$  contains the line
$(t (\xi_1)^\half, \ldots, t (\xi_k)^\half)$,
$t \in \R$,  the map is not proper.

Conversely, suppose that $\sum \xi_i \eta_i \neq 0$
whenever $\xi_i \geq 0$ are not all zero. 
Then $m = \min \{ |\Phi_H(z)| \}_{|z|^2=1}$ is positive.
Since $\Phi_H$ is quadratic, $|\Phi_H(z)| \geq m|z|^2$ for all $z$,
which implies that $\Phi_H$ is proper.
\end{proof}

This analysis already distinguishes the two possibilities
for non-empty moment fibers.

\begin{Lemma} \labell{alternative}
Let $(M,\omega,\Phi,U)$ be a proper Hamiltonian $T$-manifold.
For any $\alpha \in U$,
if the moment  fiber $\Phi\inv(\alpha)$ is not empty,
it consists of either 
\begin{enumerate}
\item \labell{case proper} 
a single orbit, which has a local model with a proper moment map, or
\item \labell{case nonproper} 
infinitely many orbits, each of which has a local model 
with a non-proper moment map.
\end{enumerate}
If $\alpha \in \text{interior}(\Phi(M))$, the second case occurs.
\end{Lemma}

\begin{proof}
Given $x \in M$, let $\rho\colon H  \to (S^1)^n$ be the
slice representation.
Consider the local model
$Y = T \times_H \C^n \times \h^0$
with   moment map 
\begin{equation} \labell{PhiY2}
\Phi_Y([t, z, \nu])  =  \alpha + \Phi_H(z) + \nu,
\end{equation}
where $\Phi_H$ is the moment map for $\rho$ and where $\alpha = \Phi(x)$.

By Lemma \ref{criterion not proper} and equations \eqref{PhiH again}
and \eqref{PhiY2}, the moment map $\Phi_Y$ is proper if and only if
the moment fiber $\Phi_Y\inv(\alpha)$ consists of a single orbit,
and otherwise $\Phi_Y\inv(\alpha)$ contains infinitely many orbits
near $\{ [t,0,0] \}$. 
The lemma now follows from the local normal form theorem 
and the connectedness of moment fibers.
\end{proof}

We will also use the following corollary of Lemma \ref{alternative}.

\begin{Lemma} \labell{121 near point}
Let $(M,\omega,\Phi,U)$ be a proper Hamiltonian $T$-manifold,
and let $\alpha$ be a point in $U$ whose
moment fiber $\Phi\inv(\alpha)$ contains exactly one orbit.
Then every neighborhood of $\alpha$ which is contained in $U$
contains a smaller neighborhood $V$ 
whose preimage, $\Phi\inv(V)$, is contractible.
\end{Lemma}

\begin{proof}
Let $Y = T \times_H \C^n \times \h^0$ be the corresponding model.
By Lemma \ref{alternative}, the moment map $\Phi_Y$ is proper.
By Proposition \ref{proper}, the preimage in $M$ and in $Y$
of a sufficiently small neighborhood $V$ of $\alpha$ are isomorphic.
Thus, we may work purely inside $Y$. Choose a neighborhood of $\alpha$
of the form $V = V_1 \times V_2$,
where $V_1 \subset \h^*$ and $V_2 \subset \h^0$ are convex,
and where we identify $\t^* = \h^* \times \h^0$.
Because $\Phi_H$ is homogeneous,
$\Phi_Y\inv(V)/T =  (\Phi_H\inv(V_1)/T) \times V_2$
is contractible. 
\end{proof}

We have already proved the local uniqueness theorem
(Theorem \ref{local uniqueness}) in case \ref{case proper} 
of Lemma \ref{alternative}. This is Proposition \ref{proper}.
Therefore, for the rest of the proof of the theorem
we may focus on case \ref{case nonproper} of Lemma \ref{alternative}.

So far, we have been allowing actions of any complexity,
but we now restrict to complexity one to define a useful polynomial:

\begin{Lemma} \labell{lem:P exact}
Let an $h$-dimensional compact abelian Lie group $H$ act on $\C^{h+1}$ 
as a subgroup of $(S^1)^{h+1}$ with a moment map that is not proper.
Then there exists a unique polynomial
\begin{equation} \labell{P}
P(z) = \prod_{j=0}^h z_j^{\xi_j},
\end{equation}
with $\xi_j \geq 0$ for all $j$, such that the following sequence is exact:
\begin{equation} \labell{P exact}
 1 \to H \stackrel{\rho}{\hookrightarrow} (S^1)^{h+1}
         \stackrel{P}{\to} S^1 \to 1.
\end{equation}
Moreover, $\xi_j > 0$ for all $j$ exactly if the moment map is onto.
\end{Lemma}

\begin{proof}[Proof of Lemma \ref{lem:P exact}]
Because the quotient $(S^1)^{h+1} / H$ is  a one dimensional compact
connected Lie group, there exists a homomorphism $P \colon (S^1)^{h+1} \to S^1$
such that the sequence \eqref{P exact} is exact.
Such a homomorphism must be of the form
\begin{equation}\labell{P lambda}
  P(\lambda) = \prod_j \lambda_j^{\xi_j}
\end{equation}
for some integers $\xi_0, \ldots, \xi_h$.
Let $\eta_j \in \h^*$ denote the weights for the $H$-action on $\C^{h+1}$.
Differentiating the identity $P \circ \rho =1$ from \eqref{P exact}, we get
\begin{equation} \labell{sum xj etaj =0}
\sum_j \xi_j \eta_j =0.
\end{equation}

By Lemma \ref{criterion not proper}, 
because the moment map $\Phi_H$ is not proper,
there exist non-negative numbers $\xi'_0, \ldots, \xi'_h$,
not all zero, such that $\sum \xi_j' \eta_j = 0$.
Since  $\sum \xi_j \eta_j = 0$, by a dimension count 
the vector $(\xi_j)$ must be a multiple of the vector $(\xi'_j)$.
Therefore, after possibly replacing the vector $(\xi_j)$ 
by the vector $(-\xi_j)$, all the $\xi_j$'s are non-negative.
By Lemma \ref{positive} and a similar dimension count, the $\xi_j$ are 
strictly positive exactly if the moment map is onto.
\end{proof}

\begin{Definition} \labell{def:P}
We call $P$ the \textbf{defining polynomial}
of the representation of $H$ on $\C^{h+1}$.
We will also use this name for
the map $P\colon \C^{h+1} \to \C$  given by the same formula
$ P(z) = \prod_j z_j^{\xi_j}$, its extension to the
local model $P\colon Y = T \times_H \C^{h+1} \times \h^0 \to \C$
given by $P([t,z,\nu]) = P(z)$, and the induced quotient map
$\ol{P} \colon Y/T \to \C$. 
We trust that this will not cause confusion. 
\end{Definition}

It is sometimes  convenient to
split a complexity one linear representation into
the direct sum of two representations,
one whose moment map is onto, and one which is toric.

\begin{Lemma}
\labell{splitting}
Let an $h$-dimensional compact abelian Lie group $H$ act on $\C^{h+1}$
through an inclusion $ \rho \colon H \hookrightarrow (S^1)^{h+1}.$
After a permutation of the coordinates, there exist splittings
$$  H = H' \times H'' \quad \quad \text{and} \quad \quad
    \C^{h+1} = \C^{h'+1} \times \C^{h''} ,$$
such that $H'$ acts on $\C^{h'+1}$ as a subgroup of $(S^1)^{h'+1}$
with a surjective moment map, and $H''$ acts on $\C^{h''}$ through an
isomorphism with $(S^1)^{h''}$.
\end{Lemma}

\begin{proof}[Proof of Lemma \ref{splitting}]
Consider the defining polynomial,
$P(z) = \prod z_j^{\xi_j}.$
Let $h''$ be the number of $j$'s such that $\xi_j=0$,
and let $h'=h-h''$.  
We can assume that $\xi_j > 0$ for $0 \leq j \leq h'$
and $\xi_{h'+j}=0$ for $1 \leq j \leq h''$.
Then $P$ defines a polynomial $P' \colon (S^1)^{h'+1} \to S^1$. 

Let us identify $H$ with its image in $(S^1)^{h+1}$.  Then
$$ H = \ker P  = \ker P' \times (S^1)^{h''}.$$
Let  $H' = \ker P'$ and $H'' = (S^1)^{h''}$.
By Lemma \ref{lem:P exact}, the moment map for $H'$ is onto.
\end{proof}

\section{The topology of the quotient}
\labell{sec:M/T}
In this section we describe the topology of the quotient $M/T$,
in preparation for showing that two such quotients are $\Phi$-diffeomorphic
if they have the same genus and isotropy data.

\begin{Proposition} \labell{topology of quotient}
Let $(M,\omega,\Phi,U)$ be a complexity one space.

The subset of $M/T$ consisting of the complement of those
moment fibers that contain single orbits is, topologically,
a manifold with boundary. 

The symplectic quotients $\Phi\inv(\alpha)/T$ 
that contain more than one point are, topologically,
closed connected oriented surfaces.
\end{Proposition}

\begin{proof}
The first claim follows immediately from Lemma \ref{alternative},
the local normal form theorem, and Lemma \ref{F homeo} below.

The fact that the symplectic quotients 
which contain more than one orbit are topological surfaces
follows immediately from Lemma \ref{alternative}, the local normal form
theorem, 
and Corollary \ref{P homeo} below. These surfaces are closed because the 
moment map is proper. 
They are connected by the connectedness of moment fibers.
The symplectic structure on the symplectic quotient induces an orientation
on the complement of a discrete set of points (namely, the exceptional orbits)
and hence on the symplectic quotient itself.
\end{proof}

\begin{Lemma} \labell{F homeo}
Let $T$ be a torus.
Let a closed $h$-dimensional subgroup $H \subseteq T$
act on $\C^{h+1}$ as a subgroup of $(S^1)^{h+1}$ with a non-proper moment map 
and a defining polynomial
$P(z) = \prod z_j ^{\xi_j}$.
Consider the model $Y = T \times_H \C^{h+1} \times \h^0$ and the map
$\ol{\Phi}_Y \colon Y/T \to \t^*$ induced by the moment map.
Define a map
$$F \colon Y / T \to \t^* \times \C$$
by
$$ F := (\ol{\Phi}_Y, \ol{P}).$$
Then $F$ is a homeomorphism of $Y/T$ with its image, which is
the polygonal set $(\image \Phi_Y) \times \C$.
\end{Lemma}

\begin{Corollary} \labell{P homeo}
The restriction of the defining polynomial to the symplectic quotient, 
$\ol{P}_\alpha \colon \Phi_Y\inv(\alpha)/T \to \C$,
is also a homeomorphism for all $\alpha \in \image \Phi_Y$.
\end{Corollary}

\begin{Definition} \labell{def F}
Let $Y = T \times_H \C^{h+1} \times \h^0$ be a local model
with a non-proper moment map. The map $F$ of Lemma \ref{F homeo}
is called the \textbf{trivializing homeomorphism} of  the model. 
\end{Definition}

The motivation for the name ``trivializing homeomorphism" is
that $F$ exhibits the quotient $Y/T$ as a trivial bundle,
with fiber $\C$, over a polygonal subset of $\t^*$.
Moreover, once we remove the moment fibers which contain
single orbits, the map $\ol{\Phi} \colon M/T \to U$ 
induced by the moment map exhibits $M/T$, topologically, as a surface bundle 
over $\image \Phi$.
This surface bundle plays an important role in the global
classification of complexity one spaces, which will be given
in subsequent papers.

\begin{proof}[Proof of Lemma \ref{F homeo}]

To show that $F$ is a homeomorphism, it is both necessary and 
sufficient to prove that the map
$(\Phi_H,P)\colon \C^{h+1} \to (\image \Phi_H) \times \C$
is onto and proper and that its fibers are exactly the $H$-orbits.
(This follows from the formulas for $F$ and $\Phi_Y$.)

We will begin by assuming that the moment map $\Phi_H$ is onto $\h^*$.
By Lemma \ref{lem:P exact}, this implies that the  $\xi_j$'s are positive.

Consider the commuting diagram
\begin{equation} \labell{FFbar}
   \begin{array}{ccc}
        \C^{h+1} & \stackrel{(\Phi_H,P)}{\to} & \h^* \times \C \\
  {\ss q_1} \downarrow \phantom{\ss{q_1}} & &
  \phantom{\ss{q_2}} \downarrow {\ss q_2} \\
        \R_+^{h+1} & \stackrel{\Fbar}{\to} & \h^* \times \R_+ ,
\end{array}
\end{equation}
where
$$q_1(z_0,\ldots,z_h) = (|z_0|^2,\ldots,|z_h|^2), \quad \quad
q_2(\alpha,\zeta) = (\alpha,|\zeta|^2),$$
and
$$ \Fbar(x_0,\ldots,x_h) =
(\half \sum_{j=0}^h x_j \eta_j , \prod_{j=0}^h x_j^{\xi_j}).$$

Let $W$ be the boundary of the positive orthant $\R_+^{h+1}$.
Since $\xi_j >0$ for all $j$, the map $(x,t) \mapsto x + t\xi$
is a homeomorphism which identifies the product $W \times \R_+$
with the orthant $\R_+^{h+1}$.
The map $\Fbar$ then becomes a map
from $W \times \R_+$  to $\h^* \times \R_+$,
given by the formula
$$ (x,t) \mapsto \left( \half\sum_{j=0}^h x_j \eta_j ,
                  \prod_{j=0}^h (x_j + t\xi_j)^{\xi_j} \right),$$
where in the first coordinate we used the equality
$\sum \eta_i(x_i+t\xi_i) = \sum \eta_i x_i $.

The map $\Fbar$ is one-to-one and onto, because
the function $x \mapsto \sum_{j=0}^h x_j \eta_j$ is a homeomorphism
from $W$ onto $\h^*$, and because for each $x \in W$, the function
\begin{equation} \labell{t}
 t \mapsto \prod_{j=0}^h (x_j + t\xi_j)^{\xi_j}
\end{equation}
from $\R_+$ to $\R_+$ is one to one and onto.
The function \eqref{t} approaches infinity uniformly in $x \in W$
as $t \to \infty$.  Therefore, $\Fbar$ is proper.

The properness of $(\Phi_H,P)$ follows from that of $\Fbar$ and $q_1$.

Let us now show that
$(\Phi_H,P)$ is onto $\h^* \times \C$.  Since $\Fbar$ is onto,
for any $(\alpha,\zeta) \in \h^* \times \C$
there exists $z \in \C^{h+1}$
such that $\Phi_H(z) = \alpha$ and $|P(z)|^2 = |\zeta|^2$.
Choose $b \in S^1$ so that $P(z) = b \zeta$.
Since the map $P \colon (S^1)^{h+1} \to S^1$ is onto,
there exists $a \in (S^1)^{h+1}$ such that $P(a) = b\inv$.
Then $(\Phi_H,P)(az) = (\alpha,\zeta)$.

Let us now show that the level sets of $(\Phi_H,P)$ are the orbits of $H$.
Suppose that $\Phi_H(z) = \Phi_H(z')$ and $P(z) = P(z')$
for some $z$ and $z'$ in $\C^{h+1}$.
Since $\Fbar$ is one to one, there exists $\lambda \in (S^1)^{h+1}$
such that $z' = \lambda z$.  We must show that $\lambda$ can be chosen
to be in $H$.
If all the coordinates of $z$ are non-zero, $P(\lambda z) = P(z)$
implies that $P(\lambda)=1$, which further implies that $\lambda \in H$,
by \eqref{P exact}.

If one of the coordinates of $z$, say $z_0$, is zero, then it is
enough to show that the $(S^1)^{h+1}$-orbit of $z$ coincides with the
$H$-orbit of $z$.  By a dimension count, it is enough to show that
the $(S^1)^{h+1}$-stabilizer of $z$ is not contained in $H$.  
Because $z_0 = 0$, the
$(S^1)^{h+1}$-stabilizer of $z$ contains the circle $(a_0, 1, \ldots,
1)$. Since $\xi_0 \neq 0$, the polynomial $P$ is not constant on this
circle. By exactness of \eqref{P exact}, this circle is not contained
in $H$.

For the general case, we may let 
$$\C^{h+1} = \C^{h'+1} \times \C^{h''} 
  \quad \text{and} \quad H = H' \times H''$$
be the  splitting into a surjective part and a toric part,
as described in  Lemma \ref{splitting}.
Then $\Phi_H (z,w) = (\Phi_{H'}(z) , \Phi_{H''} (w)).$
The map $z \mapsto (\Phi_{H'}(z),P(z))$  is proper, its
fibers are the $H'$ orbits, and it is onto $(\h')^* \times \C$, 
as we have shown above.
The map  $w \mapsto \Phi_{H''}(w)$ is a moment map for a toric action,
so it is proper and its level sets are $H''$ orbits.
Thus, the map
$$(\Phi_H,P) \colon (z,w) \mapsto (\Phi_{H'}(z) , \Phi_{H''}(w) , P(z))$$
is proper, onto $(\image \Phi_H) \times \C$,
and its level sets are the $H$-orbits.
This is precisely what we needed in order to deduce that $F$ is
a homeomorphism.
\end{proof}


\section{The smooth structure on the quotient.}
\labell{sec:smooth M/T}

In the previous section we showed that, if no moment fiber contains exactly 
one orbit, both the ordinary quotient $M/T$ and the symplectic quotient 
$\Phi\inv(\alpha)/T$ are topologically manifolds (with boundary).  
In this section we will show that they are smooth manifolds (with corners)
--- outside the exceptional orbits.

We have already given $M/T$ the quotient differentiable structure,
by specifying the sheaf of smooth functions to be the sheaf of functions 
that pull back to smooth $T$-invariant functions on $M$ 
(as in section \ref{sec:quotient}).
However, there is another way to give $M/T$ a differentiable structure:
we can cover $M$ with local models, $Y$, and take
the trivializing homeomorphisms $F \colon Y/T \to (\image \Phi_Y) \times \C$ 
(see section \ref{sec:M/T}) as local charts.  
Outside the set of exceptional orbits,
the two strategies give the same well-defined smooth structure:

\begin{Lemma} \labell{F smooth}
Let $T$ be a torus.  Let a closed $h$-dimensional subgroup $H$ of $T$
act on $\C^{h+1}$ as a subgroup of $(S^1)^{h+1}$ with a non-proper moment map.
Consider the model $Y = T \times_H \C^{h+1} \times \h^0$.
Let $E \subset Y$ be the union of the exceptional orbits.  Let 
$$F = (\ol{\Phi}_Y,\ol{P}) \colon Y / T \to  \t^* \times \C $$ 
be the trivializing homeomorphism (see Definition \ref{def F}).

Then the restriction of $F$ to $(Y \ssminus E)/T$ pulls back the sheaf
of smooth functions on $ \t^* \times \C$
onto the sheaf of smooth functions on the quotient.
\end{Lemma}

A similar statement holds for the symplectic quotients:

\begin{Corollary} \labell{P smooth}
The restriction of the defining polynomial 
$$\ol{P}_\alpha \colon (\Phi_Y\inv(\alpha) \cap (Y \ssminus E) ) /T \to \C$$
is a diffeomorphism onto its image.
\end{Corollary}

In contrast, at the exceptional orbits the trivializing homeomorphism 
and defining polynomial need not respect the differentiable structure. 

Our proof of Lemma \ref{F smooth} will use the following criterion
for non-exceptional orbits:

\begin{Lemma} \labell{exceptional}
Let an $h$ dimensional compact abelian Lie group $H$ act on $\C^{h+1}$ 
as a subgroup of $(S^1)^{h+1}$ with a surjective  moment map.  
Let $P(z) = \prod {z_j}^{\xi_j}$ be the defining polynomial.
The orbit of $z \in \C^{h+1}$ is exceptional unless 
\begin{enumerate} 
\item
$z_j \neq 0$ for all $j$, or 
\item 
there exists an index $i$ such that  $\xi_i =1$ and $z_j \neq 0$ 
for all $j \neq i$.
\end{enumerate}
\end{Lemma}

\begin{proof}
Since the moment map is onto, by Corollary \ref{nonex=free}
the $H$-orbit of $z$ is non-exceptional if and only if the stabilizer
of $z$ in $H$ is trivial.
We identify $H$ with the subgroup of $(S^1)^{h+1}$ by which it acts.
The stabilizer of $z$ then consists of those elements
$\lambda \in (S^1)^{h+1}$ such that 
$\lambda_j=1$ whenever $z_j \neq 0$
and such that 
$P(\lambda) = \prod \lambda_j^{\xi_j} = 1.$
\end{proof}

Before proving Lemma \ref{F smooth} in its full generality, we 
prove the following variant, which in particular implies the lemma
for the case that the moment map is onto.

\begin{Lemma} \labell{variant}
Let a compact abelian group $H$ act on $\C^{h+1}$ as a codimension one
subgroup of $(S^1)^{h+1}$ with a surjective moment map $\Phi_H$.
Denote by $U$ the union of the non-exceptional orbits in $\C^{h+1}$. 
For any manifold $N$, the map 
\begin{equation} \labell{F+parameters}
(\id,\ol{\Phi}_H,\ol{P}) \colon N \times (U/H) 
   \to N \times \h^* \times \C,
\end{equation}
given by
$$(n,[z]) \mapsto (n,\Phi_H(z),P(z)),$$
is a diffeomorphism with its image,
i.e., it pulls back the sheaf of smooth functions on $N \times \h^* \times \C$
onto the sheaf of smooth functions on the quotient. 
\end{Lemma}

\begin{proof}
Since $H$ acts freely on $U$, the quotient $U/H$ is naturally a smooth 
manifold (with the quotient differential structure).
The map \eqref{F+parameters} is smooth, and, by Lemma \ref{F homeo}, 
it is a homeomorphism onto its image.
Since a smooth homeomorphism between two smooth manifolds 
of the same dimension is a diffeomorphism 
exactly at those points where it is a submersion, we need to show 
that $(\d \Phi_H,\d P)|_z$ is onto for all non-exceptional $z \in \C^{h+1}$.

To show that $(d\Phi_H,dP)|_z$ is onto, it is enough to find
$\zeta \in T_z\C^n = \C^n$ such that 
$d\Phi_H|_z(\zeta) = d\Phi_H(\sqrt{-1}\zeta) = 0$ 
and $dP|_z(\zeta) \neq 0$. To see this note that, since $H$ acts freely, 
$\d \Phi_H|_z$ is onto $\h^*$.  Additionally, since $P$ is holomorphic, 
$dP|_z(\zeta)$ and $dP|_z(\sqrt{-1}\zeta) = \sqrt{-1} dP|_z(\zeta)$
form a real basis to $\C$. 

Recall that $P(z) = \prod z_j^{\xi_j}$ and
$\Phi_H(z) = \half \sum_j \eta_j z_j \zbar_j$.
Hence
$$ \d \Phi_H|_z (\zeta) 
  = \half \sum_j \eta_j (z_j \ol{\zeta}_j + \zeta_j \zbar_j) .$$

\subsubsection*{Subcase A: all the coordinates of $z$ are non-zero.}
In this case, 
$$ \d P|_z (\zeta)  = P(z) \sum_j \frac{\xi_j}{z_j} \zeta_j.$$
Let $(\zeta_j) = (\frac{\xi_j}{\zbar_j})$.  Then 
$$ \begin{array}{lcl}
 \d \Phi_H|_z(\zeta) & = & \half \sum_j \eta_j 
 \left( z_j \frac{\xi_j}{z_j} + \frac{\xi_j}{\zbar_j} \zbar_j \right) \\
 & = & \half \sum_j \eta_j (\xi_j + \xi_j) \quad = 0 \end{array}$$
by \eqref{sum xj etaj =0}, and
$$ \begin{array}{lcl}
 \d \Phi_H(\sqrt{-1}\zeta) & = & \half \sum_j \eta_j 
 \left( z_j (-\i \frac{\xi_j}{z_j}) 
        + \i \frac{\xi_j}{\zbar_j} \zbar_j \right) \\
 & = & \half \sum_j \eta_j ( -\i \xi_j + \i \xi_j ) \quad = 0,
\end{array} $$
whereas
$$ \d P|_z (\zeta) = P(z) \sum_j \frac{\xi_j }{z_j} \frac{\xi_j}{\zbar_j} 
   \neq 0.$$

\subsubsection*{Subcase B: 
one of the coordinates of $z$, say, $z_1$, is zero,
		   $\xi_1=1$, and $z_j \neq 0$ for all $j\neq 1$.}

In this case,
$$\d P|_z(\zeta) = (\prod_{j \neq 1} z_j^{\xi_j}) \zeta_1 .$$
Let $\zeta_1 = 1$ and $\zeta_j = 0$ for all $j \neq 1$.  Then
$$ \d \Phi_H|_z(\zeta) = \half \eta_1 (z_1 + \zbar_1) = 0,$$ 
and
$$ \d \Phi_H|_z(\sqrt{-1}\zeta) 
   = \half \eta_1 ( -\i z_1 + \i \zbar_1) = 0,$$ 
whereas 
$$ \d P|_z (\zeta) = (\prod_{j \neq 1} z_j^{\xi_j}) \neq 0.$$
\end{proof}

\begin{proof}[Proof of Lemma \ref{F smooth}]

Let $\C^{h+1} = \C^{h'+1} \times \C^{h''}$ and $H = H' \times H''$
be the splitting into a surjective part and a toric part, as described
in Lemma \ref{splitting}.  With this splitting, the local model is
$$ Y = T \times_{H'} \C^{h'+1} \times_{H''} \C^{h''} \times \h^0,$$
and its quotient is
$$ Y/T = (\C^{h'+1} / H') \times (\C^{h''}/H'') \times \h^0.$$
The union of the non-exceptional orbits in this quotient is
\begin{equation} \labell{nonex quotient}
 (U' / H') \times (\C^{h''}/H'') \times \h^0,
\end{equation}
where $U'$ is the union of the free orbits in $\C^{h'+1}$.
Under the identification $\t^* = (\h')^* \times (\h'')^* \times \h^0$,
the trivializing homeomorphism $F$ on \eqref{nonex quotient} is
$$ F ([z'],[z''],\nu) = ( \Phi_{H'}(z') , \Phi_{H''}(z''), \nu, P(z')), $$
where $P$ is the defining polynomial.

Lemma \ref{variant} implies that the map
$$ ([z'],[z''],\nu) \mapsto (\Phi_{H'}(z'), [z''], \nu, P(z')) $$
pulls back the sheaf of smooth functions on
$(\h')^* \times (\C^{h''}/H'') \times \h^0 \times \C$
onto the sheaf of smooth functions on
$(U'/H') \times (\C^{h''}/H'') \times \h^0$.
Therefore, it is enough to show that the map 
\begin{equation} \labell{the map}
 (\alpha,[z''],\nu,\zeta) \mapsto (\alpha,\Phi_{H''}(z''),\nu,\zeta)
\end{equation}
pulls back the sheaf of smooth functions on
$(\h')^* \times (\h'')^* \times \h^0 \times \C$
onto the sheaf of smooth functions on
$(\h')^* \times (\C^{h''}/H'') \times \h^0 \times \C.$

By a theorem of Schwartz \cite{schwarz:Topology},
any invariant smooth function 
can be expressed as a smooth function of real invariant polynomials.  
Since $H''$ acts on $\C^{h''}$ through an isomorphism with $(S^1)^{h''}$,
the ring of $H''$-invariant polynomials in $(\alpha,z'',\nu,\zeta)$
is generated by the coordinates of $\alpha$ and $\nu$,
the real and imaginary parts of $\zeta$, and $|z_1|^2, \ldots, |z_{h''}|^2$.
Finally, note that 
$$\Phi_{H''}(z_1,\ldots,z_{h''}) = A(|z_1|^2,\ldots,|z_{h''}|^2),$$
where $A \colon \R^{h''} \to (\h'')^*$ is the linear isomorphism
dual to the map $H'' \to (S^1)^{h''}$. Hence, every smooth invariant
function is the pullback via \eqref{the map} of a smooth function.
\end{proof}

\section{The associated surface}
\labell{sec:surface}

In this section we associate to a moment  fiber $\Phi\inv(\alpha)$
in a complexity one space a smooth \textbf{marked surface} whose 
underlying topological space is the symplectic quotient $\Phi\inv(\alpha)/T$.

We do not, however, define a functor from a category of complexity
one spaces to a category of marked surfaces.
For one thing, our construction depends on a choice.  More seriously,
a smooth map between two complexity one spaces does not induce a smooth
map between the associated surfaces, nor visa versa.  This is unfortunate,
because we will easily obtain an isomorphism between the marked surfaces
if the genus and isotropy data are the same.  In the next few sections
we will show that, in spite of this lack of functoriality, we can 
obtain an isomorphism of complexity one spaces from an isomorphism 
of marked surfaces. 

As we saw in the previous section, on the complement of the exceptional orbits
the symplectic quotient is naturally a smooth surface.   
Unfortunately, it is not naturally smooth near the exceptional orbits.
Nevertheless, we can use the defining polynomials of section \ref{sec:repr}
to give it a smooth structure.
However, in order to do this, we must make arbitrary 
choices: we must identify open subsets of the manifold with open subsets
of the local models.  We call these choices grommets:

\begin{Definition} \labell{def:grommet M}
Let $(M,\omega,\Phi,U)$ be a complexity one space.
A \textbf{grommet} is a $\PhiT$--diffeomorphism $\psi \colon D \to M$ 
from an open subset $D$ of a local model $Y = T \times_H \C^n \times \h^0$
onto an open subset of $M$.
\end{Definition}

The name ``grommet" is better motivated in the following definition: 

\begin{Definition} \labell{def:grommet Sigma}
Let $\Sigma$ be a smooth oriented two dimensional manifold.
A \textbf{grommet} at a point $q \in \Sigma$ is a diffeomorphism
$\varphi\colon B \to \Sigma$ from a neighborhood $B$ of the origin
in $\C$ onto an open subset of $\Sigma$, such that $\varphi$ sends 
the origin $0$ to the point $q$.
\end{Definition}

\begin{Remark} \labell{exist varphi}
Let $(M,\omega,\Phi,U)$ be a complexity one space and
consider $\alpha \in U$ such that the moment fiber $\Phi\inv(\alpha)$
contains more than one orbit.  Any grommet $\psi \colon D \to M$
with $\psi([t,0,0]) = \calO \in \Phi\inv(\alpha)$ 
induces a ``coordinate chart'' 
on the symplectic quotient $\Phi\inv(\alpha)/T$, 
that is, a homeomorphism $\varphi$ from a subset $B \subset \C$ into 
$\Phi\inv(\alpha)/T$, such that $\varphi(0) = \calO$.
Explicitly, the map
$\ol{P}_\alpha\colon (D \cap \Phi_Y\inv(\alpha))/T \to \C$
given by the defining polynomial is a homeomorphism onto its image $B$, and 
$\varphi := \ol{\psi} \circ \ol{P}_\alpha\inv\colon B \to \Phi\inv(\alpha)/T$ 
is a homeomorphism onto its image, where $\ol{\psi}$ is induced from $\psi$.
\end{Remark}

\begin{Definition} \labell{def Sigma}
Let $(M, \omega,\Phi,U)$ be a complexity one space and
$\alpha \in U$ a point whose moment fiber contains more than one orbit.
For each exceptional orbit $E_j$ in $\Phi\inv(\alpha)$, 
let $\psi_j \colon D_j \to M$ be a grommet 
such that $\psi_j([t,0,0]) = E_j$.
The \textbf{associated marked surface} consists of the following data:
\begin{enumerate}
\item
The connected oriented two-dimensional topological manifold 
$\Sigma = M_\red = \Phi\inv(\alpha)/T$. 
\item
The set of marked points $\{q_j\}$ in $\Sigma$ that 
corresponds to the set of exceptional orbits $\{ E_j \}$
in $\Phi\inv(\alpha)$.
\item
The smooth manifold structure on $\Sigma$ 
that is given by the following coordinate charts.
For each exceptional orbit $E_j$ in $\Phi\inv(\alpha)$, 
take the given grommet. 
For each non-exceptional orbit $\calO$ in $\Phi\inv(\alpha)$, 
choose an arbitrary grommet with $\psi([t,0,0]) = \calO$.
For each grommet, take the induced coordinate chart on $\Sigma$ 
as described in Remark \ref{exist varphi}.
\item
At each marked point $q_j$, the grommet on $\Sigma$ 
that is given by the coordinate chart of item 3.
\item
For each marked point $q_j$, a label consisting of the isotropy
representation at the corresponding exceptional orbit.
\end{enumerate}
\end{Definition}

The image of the grommet $\psi_j$ does not contain any of the other
exceptional orbits, $E_i$, $i\neq j$; this follows from the fact that
the model contains at most one exceptional orbit in each moment fiber
(see Lemma \ref{exceptional}).
The fact that the charts in item 3 give a well defined smooth
structure on $M/T$ follows from this and from the fact that
the smooth structures coincide on this complement 
(see Corollary \ref{P smooth}).

\section{Flattening the quotient}
\labell{sec:triv}

We are now ready to show that, after possibly replacing $M$ 
by the preimage of a small subset of $\t^*$, the quotient $M/T$ is
determined by the associated marked surface $\Sigma$.
Topologically, there is a homeomorphism from
the quotient $M/T$ to the product
$\Sigma \times (\image \Phi)$. 
Moreover, this homeomorphism
can be chosen to 
respect the smooth structure in a certain sense.
Such a homeomorphism is called a flattening; a  precise definition
is given below.

We begin by the definition for a local model.
Let $T$ be a torus, and let a closed subgroup $H \subseteq T$ act on $\C^n$ 
as a codimension one subgroup of $(S^1)^n$ with a non-proper moment map.
Consider the model
$$ Y = T \times_{H} \C^n \times \h^0,$$
with moment map $\Phi_Y \colon Y \to \t^*$.
Recall from Corollary \ref{P homeo} that the defining polynomial
$\ol{P} \colon Y/T \to \C$ restricts to a homeomorphism 
$\ol{P}_\alpha\colon \Phi_Y\inv(\alpha)/T \to \C$.

\begin{Definition} \labell{Flatten Y}
The \textbf{standard flattening} of $Y$ is the map 
$$ \delta \colon Y/T \to  (\Phi_Y\inv(\alpha)/T) \times (\image \Phi_Y)$$
given by 
$$\delta := ( (\ol{P}_\alpha \inv \circ \ol{P} ) , \ol{\Phi}_Y).$$
\end{Definition}

The standard flattening is a homeomorphism by Lemma \ref{F homeo}.

\begin{Definition} \labell{def ex sheet}
The \textbf{exceptional sheet} in the model $Y$ is the subset
$$ S := \{ [t,z,\nu] \in T \times_H \C^n \times \h^0 \ | \ P(z)=0 \}.$$
\end{Definition}
Every exceptional orbit is contained in $S$, by Lemma \ref{exceptional}.  
Thus, by Lemma \ref{F smooth} and Corollary \ref{P smooth},
the standard flattening $\delta$ is a diffeomorphism 
of $(Y \ssminus S)/T$ with its image.


Grommets were defined in Definition \ref{def:grommet M}.
Flattenings will involve grommets that are sufficiently large,
in the following sense:

\begin{Definition} \labell{wide grommet}
Let $(M,\omega,\Phi,U)$ be a complexity one space.
A grommet $\psi \colon D \to M$ is \textbf{wide} 
if $D$ is an open subset of a model $Y$ on which the moment map $\Phi_Y$
is non-proper and if the domain $D$ contains that part of 
the exceptional sheet that lies over $U$,
i.e., if $(\Phi_Y\inv(U) \cap S) \subset D$.
\end{Definition}

We are now ready to define the flattening of a complexity one space.

\begin{Definition} \labell{def:triv M}
Let $(M,\omega,\Phi,U)$ be a complexity one space.
Assume that the moment  fiber over $\alpha \in U$ contains 
more than one orbit.
A \textbf{flattening} of the space $(M,\omega,\Phi)$ about $\alpha$
consists of the following data.
\begin{enumerate}
\item 
A homeomorphism 
\begin{equation} \labell{delta}
 \delta \colon M/T \to (\Phi\inv(\alpha)/T)  \times (\image \Phi)
\end{equation}
whose second coordinate is induced by the moment map.
\item
For each exceptional orbit $E_j$ in $\Phi\inv(\alpha)$, a wide grommet 
$\psi_j \colon D_j \to M$ such that $\psi_j([t,0,0]) = E_j$.
\end{enumerate}
We require that the following two conditions be satisfied:
\begin{enumerate}
\item
The restriction of $\delta$ to the complement of the exceptional sheets, 
$$\delta\colon M/T \ssminus \sqcup_j \psi_j(S_j \cap D_j)/T \to
(\Phi\inv(\alpha) \ssminus \sqcup_j E_j)/T \times (\image \Phi),$$
must be a diffeomorphism,
in the sense discussed at the beginning of section \ref{sec:quotient}.
\item
Additionally, near the exceptional sheets $\delta$ must be given by the 
standard flattenings of the local models.
More precisely, the following diagram must commute: 
\begin{equation} \labell{triv compatibility}
\begin{array}{ccccc}
D_j/T  & \stackrel{\delta_j}{\to}  & 
(\Phi_j\inv(\alpha) \cap D_j)/T \times \t^* \\
 \phantom{\ss \ol{\psi}_j} \downarrow {\ss \ol{\psi}_j}
 & &
 \phantom{\ss (\ol{\psi}_j, \id,} \downarrow {\ss (\ol{\psi}_j,\id)}\\
 M/T & \stackrel{\delta}{\to} & (\Phi\inv(\alpha)/T) \times \t^*, \\
\end{array}
\end{equation}
where $\Phi_j$ is the moment map on the corresponding 
local model $Y_j \supset D_j$, where  $\delta_j \colon Y_j/T  \to 
   (\Phi_j\inv(\alpha)/T)  \times (\image \Phi_j)$
denotes the standard flattening of $Y_j$,
and where $\ol{\psi}_j \colon D_j/T \to M/T$ is induced by the grommet.
In particular, we require that the diagram  be well defined,
i.e., that the image $\delta_j(D_j/T)$ be contained in
$((\Phi_j\inv(\alpha) \cap D_j)/T) \times \t^*$.
\end{enumerate}
\end{Definition}

The rest of this section is devoted to showing that flattenings 
always exist locally:

\begin{Lemma} \labell{exists local triv}
Let $(M,\omega,\Phi,U)$ be a complexity one space,
and let $\alpha \in \t^*$ be a point 
whose moment fiber contains more than one orbit. 
Then there exists a neighborhood $V$ of $\alpha$ contained in $U$
whose preimage, $\Phi\inv(V)$, admits a flattening about $\alpha$.
\end{Lemma}

Since $\image \Phi$ is convex, and hence contractible,
the following consequence is immediate:

\begin{Corollary} \labell{cor 121}
Let $(M,\omega,\Phi,U)$ be a complexity one space.
Consider $\alpha \in U$ such that the moment fiber
$\Phi\inv(\alpha)$ contains more than one orbit.  Then 
for every sufficiently small convex neighborhood $V$ of $\alpha$,
the restriction map
$$ H^*(V/T) \to H^*(\Phi\inv(y)/T)  $$
is an isomorphism for all $y \in \Phi(V)$. 
In particular, $\Phi\inv(V)$
satisfies Condition \eqref{technical121}.
\end{Corollary}    

Since the set of points in $U$ whose moment fiber has more than one
orbit is connected by Lemma \ref{alternative}, Lemma \ref{exists local triv}
has the following additional immediate consequence: 

\begin{Corollary} \labell{same genus}
Let $(M,\omega,\Phi,U)$ be a complexity one space.
Then all the symplectic quotients $\Phi\inv(\alpha)/T$  
that contain more than one point have the same genus.
Thus, the genus of a complexity one space 
(see section \ref{sec:intro}) is well-defined.
\end{Corollary}

The following lemma is a first step for the proof 
of Lemma \ref{exists local triv}.

\begin{Lemma} \labell{exist wide grommets}
Let $(M,\omega,\Phi,U)$ be a complexity one space,
and assume the moment fiber over $\alpha \in U$ contains 
more than one orbit.  Denote the exceptional orbits in
$\Phi\inv(\alpha)$ by $\{ E_j \}$. 

After replacing $M$ by the preimage of some neighborhood of $\alpha$ 
in $U$, there exist wide grommets $\psi_j  \colon D_j \to M$ 
such that $\psi_j([t,0,0]) = E_j$ and the images $\psi_j(D_j)$ 
have pairwise disjoint closures.
\end{Lemma}

\begin{proof}
By Lemma \ref{alternative}, for every exceptional orbit $E_j$
over $\alpha$, the corresponding local model has a non-proper
moment map, $\Phi_j  \colon Y_j \to \t^*$.
By the local normal form theorem, we may choose 
a grommet $\psi_j  \colon D_j \to M$ such that $\psi_j([t,0,0]) = E_j$.

By Lemma \ref{F homeo} and Definition \ref{def ex sheet}, 
the moment map $\Phi_j$ restricts to a 
homeomorphism of the exceptional sheet  
$S_j \subset Y_j$ with the image of $\Phi_j$.
Hence there exists a neighborhood $W_j$ of $\alpha$ such that 
$S_j \cap D_j = S_j \cap \Phi_j\inv(W_j)$.

For $i \neq j$, the intersection 
$\psi_i(S_i \cap D_i) \cap \psi_j(S_j \cap D_j)$
is a closed subset of $M$ which does not meet the fiber $\Phi\inv(\alpha)$. 
Since the moment map is proper, there exists a neighborhood $V$ 
of $\alpha$ which does not meet the image 
under the moment map of any of these intersections. 

If we define  $W := \bigcap\limits_j W_j \cap V$ of $\alpha$, and 
replace $M$ by $M \cap \Phi\inv(W)$ and $D_j$ by $D_j \cap \Phi_j\inv(W)$, 
the grommets $\psi_j$ become wide.
Also, the exceptional sheets $\psi_j(S_j \cap D_j)$ are then closed 
and disjoint, so we can shrink each $D_j$ to a smaller neighborhood 
of $S_j \cap D_j$ to obtain wide grommets whose images have pairwise 
disjoint closures.
\end{proof}

\begin{proof}[Proof of Lemma \ref{exists local triv}]
Let $\psi_j  \colon D_j \to M$ be wide grommets
such that $\psi_j([t,0,0]) = E_j$ are the exceptional orbits
in the moment fiber $\Phi\inv(\alpha)$
and such that the images $\psi_j(D_j)$ have disjoint closures in $M$.
These grommets exist by Lemma \ref{exist wide grommets}.
This will not be ruined if we further restrict to a smaller neighborhood 
of $\alpha$.

Recall that the standard flattening of the local model $Y_j$ is
$$ \delta_j = (g_j , \ol{\Phi}_j)  \colon Y_j / T 
   \to (\Phi_j\inv(\alpha) / T) \times (\image \Phi_j) $$
where $g_j = (\ol{P}_{j,\alpha})\inv \circ \ol{P}_j$.
Replace  $D_j / T$ by its intersection with 
$ g_j\inv ( (\Phi_j\inv(\alpha) \cap D_j) / T ) .$
Then the restriction
$$	\delta_j  \colon D_j / T 
	\to ( (\Phi_j\inv(\alpha) \cap D_j) / T ) \times ( \image \Phi_j) $$
is well defined. After this, the grommets determine a unique map $\delta$
on the images of $D_j/T$ in $M/T$ such that the following diagram commutes.

\begin{equation} \labell{triv compatibility1}
\begin{array}{ccccc}
 D_j/T  & \stackrel{\delta_j}{\to} &  
 ((\Phi_j\inv(\alpha) \cap D_j)/T) \times \t^* \\
 \phantom{\ss \ol{\psi}_j} \downarrow {\ss \ol{\psi}_j}
 & &
 \phantom{\ss (\ol{\psi}_j, \id )} \downarrow {\ss (\ol{\psi}_j,\id)}\\
 \bigsqcup_j \ol{\psi}_j (D_j/T) & 
\stackrel{\delta}{\to} & (\Phi\inv(\alpha)/T)  \times \t^*.
\end{array}
\end{equation}
We need to extend $\delta$ to the rest of $M/T$, 
perhaps after shrinking the $D_j$s 
to smaller neighborhoods of $S_j \cap \Phi_j\inv(U)$.

Using the stability of the moment map, Lemma \ref{F smooth} implies that
on the complement of the exceptional sheets in the quotient $M/T$,
the map
$\ol{\Phi} \colon 
 M/T \ssminus \sqcup_j \psi_j(S_j \cap D_j)/T \to (\image \Phi)$ 
induced by the moment map is a submersion.
Namely, for each point $[m]$ in the domain of this map 
there exists a neighborhood $W$ of $\Phi(m)$ in $\t^*$
such that a neighborhood of $[m]$ is diffeomorphic to 
the product of a disk with $W \cap (\image \Phi)$
with the map $\ol{\Phi}$ being the projection map.

The partial flattening \eqref{triv compatibility1} determines
an Ehresmann connection for this submersion, defined on the open subset
$\sqcup_j \psi_j (D_j \ssminus S_j)/T$:  
we declare the horizontal tangent vectors to be those whose push-forward 
by $\delta$ is tangent to the sheets 
$\{ q \} \times \t^*$ for $q \in \Phi\inv(\alpha)/T$.

We extend this to an Ehresmann connection on the entire complement of the 
exceptional sheets, $M/T \ssminus \sqcup_j \psi_j(S_j\cap D_j)/T$, 
perhaps after shrinking the $D_j$s; this is easily done with a partition 
of unity. Then for a point $p$ in 
$M/T \ssminus \sqcup_j \psi_j (S_j \cap D_j) /T $,
any path $\gamma$ in $U$ which starts at $\ol{\Phi}(p)$ can 
be lifted to a horizontal path in 
$M/T \ssminus \sqcup_j \psi_j (S_j \cap D_j)/T $.

We proceed as in the proof of Ehresmann's lemma.
Let us assume that $\alpha=0$ and that $U$ is a ball centered at $0$.
We can choose coordinates on $\t^*$ such that $(\image \Phi)$ becomes
$$ (\image \Phi) = U \cap (\R^k \times \R_+^l) \ , \ k + l = m = \dim \t^*.$$
(This is possible by Lemma \ref{nonexceptional near boundary}.)
Denote by $v_1, \ldots, v_m$ the standard vector fields on $\t^*$
that are parallel to the coordinate axes, let 
$\tilde{v}_1, \ldots, \tilde{v}_m$ be their horizontal liftings to
$M/T \ssminus \sqcup_j \psi_j (S_j \cap D_j)/T$,
and let $f_j^t$, for $i=1, \ldots, m$ and $t \in \R$, 
be the flows which the $\tilde{v}_i$ generate.  For $p \in M/T$, define
$$	\delta(p) = (g(p) , \ol{\Phi}(p)) $$
where, if $(t_1,\ldots,t_m)$ are the coordinates of $p$,
then $g(p) \in \Phi\inv(\alpha)/T$ is given by
$g(p) = f_1^{-t_1} \cdots f_m^{-t_m} (p)$.
\end{proof}

\section{Diffeomorphism between quotients}

In this section, we show that if two complexity one
spaces equipped with flattenings
have the same genus and isotropy data, their quotients
are $\Phi$-diffeomorphic, in the sense of Definition \ref{Phi diffeo}:

\begin{Proposition} \labell{rigid extends}
Let $(M,\omega,\Phi,U)$ and $(M',\omega',\Phi',U)$ be 
complexity one spaces that admit flattenings about a point $\alpha \in U$. 
If the spaces have the same genus and isotropy data,
then there exists a $\Phi$-diffeomorphism from $M/T$ to $M'/T$.
\end{Proposition}

\begin{Definition}
Let $\Sigma$ and $\Sigma'$ be closed oriented surfaces
with labeled marked points and with grommets at these points.
(See Definition \ref{def:grommet Sigma}.)
An orientation preserving diffeomorphism $g \colon \Sigma \to \Sigma'$ 
is \textbf{rigid} if
\begin{itemize}
\item
it induces a bijection between the marked points in $\Sigma$ 
and those in $\Sigma'$, and sends each marked point $q_j$ in $\Sigma$
to a marked point $q'_j$ in $\Sigma'$ with the same label;
\item
for each marked point $q_j \in \Sigma$ and $q'_j \in \Sigma'$
and corresponding grommets $\varphi_j$ and $\varphi'_j$, the composition
$${\varphi'_j}\inv \circ g \circ \varphi_j,$$ 
which is a diffeomorphism from a neighborhood of $0$ in $\C$
to another neighborhood of $0$ in $\C$, coincides with a rotation of $\C$
on some (smaller) neighborhood of $0$.
\end{itemize}
\end{Definition}

\begin{proof}[Proof of Proposition \ref{rigid extends}]
Let $(M,\omega,\Phi,U)$ and $(M',\omega',\Phi',U)$ be 
complexity one spaces with flattenings about the point $\alpha \in U$. 
Let $\Sigma$ and $\Sigma'$ be the associated marked surfaces,
as in Definition \ref{def Sigma}.

If the spaces have the same isotropy data, 
there exists a bijection from the marked points in $\Sigma$ 
onto the marked points in $\Sigma'$ which respects the isotropy labels.
If the spaces have the same genus, this bijection extends
to an orientation preserving diffeomorphism from $\Sigma$ to $\Sigma'$;
this follows from standard differential topology.
Moreover, this diffeomorphism can be deformed near the marked points
into a rigid diffeomorphism $h \colon \Sigma \to \Sigma'$.
This type of result is standard in differential topology; see, e.g., 
\cite[II, 5.2]{kosinski}. 

The images of the moment maps $\Phi$ and $\Phi'$ are the same;
this follows from Lemma \ref{nonexceptional near boundary}
and from the fact that the isotropy data are the same.
Let
$$ \delta  \colon M/T \to \Phi\inv(\alpha)/T  \times (\image \Phi) $$
and
$$\delta'  \colon M'/T \to {\Phi'}\inv(\alpha)/T  \times (\image \Phi') $$
be maps given in the flattenings. 
Under the identification of the symplectic quotients
$\Phi\inv(\alpha)/T$ and ${\Phi'}\inv(\alpha)/T$
with $\Sigma$ and $\Sigma'$, respectively, 
the rigid diffeomorphism $h \colon \Sigma \to \Sigma'$ extends to a map
$$ (h, \id)  \colon 
   (\Phi\inv(\alpha)/T) \times (\image \Phi)
   \to ({\Phi'}\inv(\alpha)/T) \times (\image \Phi') .$$
We will show that the map $g \colon M/T \to M'/T$ defined by 
$$g := {\delta'}\inv \circ (h,\id) \circ \delta$$
is a $\Phi$-diffeomorphism.

The diffeomorphism $h$ fixes an identification
between exceptional orbits in $\Phi\inv(\alpha)$ and ${\Phi'}\inv(\alpha)$
with the same isotropy representation.
Thus, we can unequivocally denote by $\{Y_j\}$ the local models for the 
exceptional orbits over $\alpha$ in both $M$ and $M'$.  Let 
$$\psi_j  \colon D_j \to M  
  \quad \text{and} \quad  \psi'_j  \colon  D'_j \to M' $$
denote the grommets,
with $D_j \subseteq Y_j$ and $D'_j \subseteq Y_j$,
and let $E_j$ and $E_j'$ denote the exceptional orbits
in $\Phi\inv(\alpha)$ and in ${\Phi'}\inv(\alpha)$.

Our first claim is that the restriction
$$ g \colon \left( M \ssminus \sqcup_j \psi_j(S_j \cap D_j) \right)/T 
   \to \left(M' \ssminus \sqcup_j \psi_j'(S_j \cap D_j) \right)/T$$
is a $\Phi$-diffeomorphism.
This is easy: by the definition of flattening, the restrictions 
$$
\delta \colon \left( M \ssminus \sqcup_j \psi_j(S_j \cap D_j) \right) /T 
\to 
\left( \Phi\inv(\alpha) \ssminus \sqcup_j E_j \right)/T 
\times (\image \Phi) $$
and
$$\delta' \colon \left( M' \ssminus \sqcup_j \psi'_j(S_j \cap D_j) \right) /T 
\to 
\left( {\Phi'}\inv(\alpha) \ssminus \sqcup_j E'_j \right)/T 
\times  (\image \Phi') $$
are both diffeomorphisms.  Moreover, the map
$$ \left( \Phi\inv(\alpha) \ssminus \sqcup_j E_j \right)/T  
\to  \left( {\Phi'}\inv(\alpha) \ssminus \sqcup_j E'_j \right)/T  $$ 
induced by $h$ is a diffeomorphism,
since the smooth structures on $\Phi\inv(\alpha)/T$ and $\Sigma$ 
agree off the exceptional orbits.

It remains to show that $g$ is a $\Phi$-diffeomorphism in
a neighborhood of each exceptional sheet $\psi_j(S_j \cap D_j)/T$.

Let $\varphi_j \colon B_j \to \Sigma$ and $\varphi'_j \colon B'_j \to \Sigma'$
denote the grommets of the associated surfaces.
Since $h$ is rigid,  there exist
$a_j \in S^1$ such that $\varphi_j\inv \circ h \circ \varphi_j$ 
is given by  rotation by $a_j \in S^1$ on some
neighborhood of the origin in $\C$ 

Let $P_j \colon (S^1)^{n_j}  \to S^1$ be the defining polynomial 
for the exceptional orbit $E_j$.  Since $P_j$ is surjective,
we may choose  $\lambda_j \in (S^1)^{n_j}$ so that $P_j(\lambda_j) = a_j$.
This defines an  equivariant symplectomorphism from the local model 
$Y_j = T \times_{H_j} \C^{n_j} \times \h_j^0$ to itself as follows:
\begin{equation} \labell{eq:flat}
 \lambda_j \cdot ([t,z,\nu]) = [ t, \lambda_j \cdot z, \nu].
\end{equation}
This map induces a $\Phi$-diffeomorphism on the quotient,
 $g_j\colon Y_j/T \to Y_j/T$.
It remains to show only that the $g_j$ and $g$ agree in
some neighborhood of $\psi_j(D_j \cap S_j)$. 
Indeed, when we use the trivializing homeomorphism $F_j$
to identify  $Y_j$ with $\C \times (\image \Phi_j)$,
the map $g_j$ sends $(z,\beta)$ to $(a_j z, \beta)$.
\end{proof}


\section{Proof of the  Local Uniqueness Theorem.}

We now have all the ingredients to prove Theorem 
\ref{local uniqueness}. We recall the statement:

\medskip \noindent
\textbf{Theorem \ref{local uniqueness}.}
\emph{
Let $(M,\Phi,\omega,U)$ and $(M',\Phi',\omega',U)$ be complexity one spaces.
Assume that their \DH measures are the same, and that their genus 
and isotropy data over a point $\alpha \in \t^*$ are the same. 
Then there exists a neighborhood of the point $\alpha$ 
over which the spaces are isomorphic.
} \smallskip

\begin{proof}
Since the case that the moment fiber  
$\Phi\inv(\alpha)$ contains just one orbit is covered
by  Proposition \ref{proper},
we may assume that the moment fiber contains more than one orbit.

By Lemma \ref{exists local triv},  
after possibly restricting to the preimage of a smaller neighborhood 
of $\alpha$, 
we may assume that $M$ and $M'$ are equipped with flattenings.
By assumption, the spaces $M$ and $M'$ have the same genus
and isotropy data.
Therefore, by Proposition \ref{rigid extends},
there is a $\Phi$-diffeomorphism $g\colon M/T \to M'/T$.

Since the spaces have flattenings, Condition 
\eqref{technical121} is satisfied. (See Corollary \ref{cor 121}.)
By assumption, the \DH measures of $M$ and $M'$ are the same.
Hence, we can apply Propositions \ref{eliminate} and \ref{prop:HS}.
The first implies that the map $g$ 
lifts to a $\PhiT$--diffeomorphism from $M$ to $M'$.
The second then guarantees that there exists 
$\PhiT$--symplectomorphism from $M$ to $M'$.
\end{proof}

\section{Proof of uniqueness for centered spaces}
\labell{sec:centered}

In this section we prove Theorem \ref{thm:centered-uniqueness}, which 
shows that the invariants that we described also separate centered spaces.

We recall Definition \ref{centered-definition}:
A proper Hamiltonian $T$-manifold $(M,\omega,\Phi,U)$ 
is \textbf{centered} about a point $\alpha \in U$ 
if $\alpha$ is contained in the closure of the moment image 
of every orbit type stratum in $M$.
The local normal form theorem, together with the properness of the moment map,
imply that every point in $\t^*$ has a neighborhood whose preimage
is centered. 

We recall the statement of the theorem.

\smallskip \noindent 
\textbf{Theorem \ref{thm:centered-uniqueness}} (Centered Uniqueness).
\emph{
Let $(M,\Phi,\omega,U)$ and $(M',\Phi',\omega',U)$ be complexity one spaces 
that are centered about $\alpha \in U$.  
Assume that their \DH measures are the same and that their genus and 
isotropy data over $\alpha \in \t^*$ are the same. Then the spaces are 
isomorphic.}
\smallskip

\begin{Remark} \labell{easy centered}
In Theorem \ref{thm:centered-uniqueness},
if the moment fibers $\Phi\inv(\alpha)$ and
${\Phi'}\inv(\alpha)$ are each  a single orbit,
the centered spaces are isomorphic  if the isotropy data 
are the same. (The \DH measures 
and the genus are then  automatically the same.)
The proof for this case only uses three results from earlier sections:
Propositions \ref{proper} and \ref{eliminate} and Lemma \ref{alternative}.
\end{Remark}

The following proof of Theorem \ref{thm:centered-uniqueness}
relies on a couple of technical lemmas which we postpone until 
after the proof.

\begin{proof}
\ \\
\noindent 
\textbf{Case I: the moment  fiber is a single orbit.}
By Proposition \ref{proper}, there exists a convex sub-neighborhood 
$V \subset U$ of $\alpha$ and a $\PhiT$--diffeomorphism 
(in fact, symplectomorphism) from $\Phi^{-1}(V)$ to ${\Phi'}^{-1}(V)$.
By Lemma \ref{121 near point} we can choose $V$ so that
$\Phi\inv(V)$ and ${\Phi'}\inv(V)$ satisfy Condition \eqref{technical121}.

By Lemma \ref{stretching}, this implies that there exists a 
$\PhiT$--diffeomorphism from $(M,\omega,\Phi)$ to $(M',\omega',\Phi')$, and 
that $M$ and $M'$ themselves also satisfy Condition \eqref{technical121}.
The \DH measures coincide;
hence we may apply Proposition \ref{eliminate}, which completes the proof.

\noindent 
\textbf{Case II: the moment fiber contains more than one orbit.}
By Lemma \ref{exists local triv},  
there exists a convex sub-neighborhood $V \subset U$ of $\alpha$  
so that $\Phi\inv(V)$ and ${\Phi'}\inv(V)$ are equipped with flattenings.
By assumption, the spaces $M$ and $M'$ have the same genus
and isotropy data.
Therefore, by Proposition \ref{rigid extends},
there is a $\Phi$-diffeomorphism $g\colon \Phi\inv(V) \to {\Phi'}\inv(V)$.
Since these spaces have flattenings, Condition \eqref{technical121} 
is satisfied. (See Corollary \ref{cor 121}.)
By assumption, the \DH measures of $M$ and $M'$ are the same.
Hence,  Proposition  \ref{prop:HS} implies that the map $g$ 
lifts to a $\PhiT$--diffeomorphism from $\Phi\inv(V)$ to ${\Phi'}\inv(V)$.
Proposition \ref{eliminate} completes the proof.
\end{proof}

We have used the following ``stretching lemma", which tells us that 
a centered space retracts onto a neighborhood of its central fiber.
In a later paper we will use this lemma for spaces which may have
complexity greater than one. 

\begin{Lemma} \labell{stretching}
Let $\t^*$ be the dual of the Lie algebra of a torus $T$,
$U \subset \t^*$ an open convex neighborhood 
of a point $\alpha \in \t^*$, and $V \subset U$ any sub-neighborhood.
Then there exists a convex neighborhood $W$ of $\alpha$ contained in $V$,
and a diffeomorphism $f \colon U \to W$
with the following property:
for any proper Hamiltonian $T$-manifold $(M,\omega,\Phi,U)$  
that is centered about $\alpha$, 
there exists a smooth equivariant orientation preserving diffeomorphism
$F\colon M \to  \Phi^{-1}(W)$ such that $\Phi \circ F  = f \circ \Phi$. 
\end{Lemma}

Before beginning the proof of Lemma \ref{stretching},
recall that the \textbf{Euler vector field} on a vector space $V$ is given by
$X = \sum x_i \dd{x_i}$, where $x_i$ are linear coordinates. 
This vector field is the generator of the flow $x \mapsto e^t x$, 
thus it is independent of the choice of coordinates.

\begin{Lemma} \labell{radial}
Let $(M,\omega,\Phi,U)$ be a proper Hamiltonian $T$-manifold.
Suppose that $U$ contains the origin $0$ of $\t^*$, 
and that the space is centered about the origin.
Then the Euler vector field $X$ on $\t^*$ lifts to a smooth invariant
vector field $\tilde{X}$ on $M$, that is, $\Phi_*(\tilde{X}) = X$.
\end{Lemma}

\begin{proof}
By the local normal form theorem, it is enough to construct 
the vector field $\tilde{X}$ on the local models.
We can then patch together the pieces by an invariant partition of unity.

Notice that if a map $\Phi \colon V \to W$ between vector spaces 
is homogeneous of degree $m$, then $\Phi_* X_V = m X_W$ where
$X_V$ and $X_W$ are the Euler vector fields on $V$ and $W$;
this follows from the equality $\Phi(e^t v) = e^{mt} \Phi(v)$. 
In particular, the Euler vector field on $W$ lifts to a vector field on $V$.
Similarly, if $\Phi_i \colon V_i \to W_i$, $i=1,2$, are homogeneous
(possibly of different degrees), and
$\Phi = \Phi_1 \times \Phi_2 \colon V_1 \times V_2 \to W_1 \times W_2$,
then the Euler vector field on $W_1 \times W_2$ lifts to a 
vector field on $V_1 \times V_2$.

Consider a local model in $M$, namely,
$Y = T \times_H \C^n \times \h^0$, with a moment map
$\Phi_Y ([t,z,\nu]) = \alpha + \Phi_H(z) + \nu$. 
The stratum fixed by $H$ is $T \times_H (\C^n)^H \times \h^0$.
Because the space is centered, we must have that $\alpha \in \h^0$.
Without loss of generality we may assume that $\alpha=0$.
To lift the Euler vector field on $\t^*$ to a $T$-invariant vector field 
on $Y$, it is enough to lift it to an $H$-invariant vector field
on $\C^n \times \h^0$.  This is possible by the previous paragraph,
because $\Phi|_{\C^n \times \h^0}$ is bihomogeneous.
\end{proof}

\begin{proof} [Proof of Lemma \ref{stretching}.] 
Without loss of generality we may assume that $\alpha = 0$.
Choose $\epsilon > 0 $ so that an $\epsilon$-ball about $0$
is contained in $V$.
Let  $g_t \colon [0,\infty) \to [0,\infty)$ for $0 \leq t \leq 1$  be
an isotopy such that $g_0$ is the identity map, 
the image of $g_1$ is contained in $[0,\epsilon)$, 
$g_t(x) \leq x$ for all $x$ and all $t$, and $g_t(x) = x$
for all $x$ near zero and all $t$.

Take $f_t(v) = g_t(|v|) \frac{v}{|v|}$ for all $v \in U \ssminus \{ 0 \}$
and $f_t(0) = 0$; let $f = f_1$.

Let $\xi_t$ be the vector field on $V$ which generates this isotopy:
$\frac{df_t}{dt} = \xi_t \circ f_t$. 
Since $\xi_t$ vanishes near $v=0$, we can
write $\xi_t = \psi_t \cdot X$, where $\psi_t\colon \t^* \to \R$ is a smooth
function, and $X$ is the Euler vector field on $\t^*$.
By Lemma \ref{radial}, there exists a smooth invariant
vector field $\tilde{X}$ on $M$ such that $\psi_*(\tilde{X}) = X$.
So  $\tilde{\xi}_t = (\psi_t \circ \Phi) \cdot \tilde{X}$ is a smooth 
invariant vector field on $M$ which is a lifting of  $\xi_t$. 
Because $\Phi$ is proper, the vector field $\tilde{\xi}_t$ generates 
an isotopy, $F_t$.  Take $F = F_1$.
\end{proof}

\section{Application to packings of Grassmanians} 
\labell{sec:applications}

We are now ready to  present our  application.
First, we recall a definition from symplectic topology:

\begin{Definition}
A symplectic manifold $M$ admits a \textbf{full packing 
by $\mathbf{k}$ equal balls} if for any $\epsilon > 0$ 
there exists a symplectic embedding into $M$
of a disjoint union of $k$ symplectic balls with equal radii 
such that the complement of the image has volume less than $\epsilon$.
\end{Definition}

Let  $\Gr^+(2,\R^n)$ denote
the Grassmanian of all oriented real 2-planes in $\R^n$,
together with its canonical (up to multiplication by a constant)
$\SO(n)$-invariant symplectic structure, and the 
$\lfloor \frac{n}{2} \rfloor$ dimensional torus action given 
by restricting the standard action of $\SO(n)$.

\medskip \noindent
\textbf{Theorem \ref{thm:applications}.}
\emph{
Let $M$ be the Grassmanian $\Gr^+(2, \R^5)$ or $\Gr^+(2, \R^6)$. 
There exists an equivariant symplectic embedding of a disjoint union of
two symplectic balls with linear actions and with equal radii into $M$
such that the complement of the image has zero volume.
A fortiori, these  Grassmanians 
can be fully packed by two equal balls.} \smallskip

The following tool will be useful:

\begin{Lemma} \labell{lem:ellipsoid}
Let $(M,\omega,\Phi)$ be a compact complexity one space over $\t^*$.
Let $p \in M$ be an isolated fixed point
with isotropy weights $\eta_1,\ldots,\eta_n$.
Assume that the differences $\eta_i - \eta_j$ span 
a codimension one subspace, $H$, of $\t^*$.
Assume, moreover, that $p$ is the only fixed point
whose moment map image lies on one open side,  $H_+$, of $H$.

Then the preimage $\Phi\inv(H_+)$ is equivariantly symplectomorphic 
to a ball with a linear $T$-action.
\end{Lemma}
\begin{proof}
First, we show that $\Phi\inv(H_+)$ is centered.
The closure $N$  of  an orbit type stratum  in $M$ 
is itself a compact symplectic manifold 
with the restricted $T$ action and moment map. By the convexity theorem, 
its moment image is the convex hull of the moment images of its fixed points. 
Either  $N$ contains $p$, or its moment image would be contained in 
$\conv(M^T \ssminus p)$, and therefore disjoint from $H_+$.

Let $T$ act on $\C^n$ with weights $\eta_1, \ldots, \eta_n$
and with the moment map that sends $(z_1,\ldots, z_n)$ to
$\Phi(p) + \sum \half \eta_i |z_i|^2$.  
The moment preimage of $H_+$ in $\C^n$ is a ball.
Hence, this ball is also a centered complexity one space over $H_+$.

Since both spaces are centered about $a=\Phi(p)$, and
the preimages of $a$ are both single orbits with the same isotropy data,
by Theorem \ref{thm:centered-uniqueness} the spaces are 
equivariantly symplectomorphic.
\end{proof}
Now consider any semi-simple compact Lie group $G$, and
let $T$ be a maximal torus.
Use the Killing form to identify $\t$ and $\t^*$ and embed $\t^*$ in $\g^*$.
Recall that the coadjoint orbit in $\g^*$ through an element $x$ of $\t^*$
is a symplectic manifold,
and the projection to $\ft^*$ is a moment map for the $T$ action.
The fixed points for the $T$-action 
are exactly the Weyl group orbit of $x$ in $\t^*$.
Moreover, the isotropy weights at a fixed point $y \in \t^*$ 
are exactly those roots $\alpha \in \t^*$ 
for which $\left< \alpha , y \right> < 0$.
\begin{proof}[Proof of Theorem \ref{thm:applications} for $\Gr^+(2,\R^5)$]
The Lie algebra of the  maximal torus of $\SO(5)$ 
can be identified with $\R^2$ with the standard metric.
The roots are $(\pm 1, 0)$, $(0,\pm 1)$, $(\pm 1, \pm 1)$.
The Weyl group acts by permuting the coordinates and by flipping
their signs. 

The orbit through the point $(1,0)$ is naturally identified 
with the Grassmanian $\Gr^+(2,\R^5) = \SO(5) / S(O(2) \times O(3))$.
The Weyl group orbit of this point consists of the points
$(1,0)$, $(-1,0)$ $(0,1)$, and $(0,-1)$.
The image of this moment map is a diamond.
The isotropy weights at $(1,0)$ are $(-1,1)$, $(-1,0)$, and $(-1,-1)$.
By Lemma \ref{lem:ellipsoid},
the preimage of the half space $\{ (x,y) \ | \ x>0 \}$ is a ball
as required. A similar argument shows that the preimage of
the opposite half space is again a ball.
\end{proof}

\begin{proof}[Proof of Theorem \ref{thm:applications} for $\Gr^+(2,\R^6)$]

The Lie algebra of the  maximal torus of $\SO(6)$ 
can be identified with $\R^3$ with the standard metric.
The Weyl group acts by permuting the coordinates and by flipping
the signs of two coordinates at a time.
The roots are $(\pm 1, \pm 1, \pm 1)$.

The orbit through the point $(1,0,0)$ is naturally identified 
with the Grassmanian $\Gr^+(2,\R^6) = \SO(6) / S(O(2) \times O(4))$.
The Weyl group orbit of this point consists of the points
$(\pm 1,0,0)$, $(0,\pm 1, 0)$, and $(0,0,\pm 1)$.
The image of this moment map is an octahedron.
The isotropy weights at $(1,0,0)$ are $(-1,\pm 1, \pm 1)$.
By Lemma \ref{lem:ellipsoid}, 
the preimage of the half space $\{ (x,y,z) \ | \ x>0 \}$ is a ball
as required. A similar argument shows that the preimage of
the opposite half space is again a ball.
\end{proof}


\begin{thebibliography}{KKMS}   

\bibitem[AH]{ah-ha} K. Ahara and A. Hattori,
\emph{ 4 dimensional symplectic $S^1$-manifolds admitting moment map}, 
J.\ Fac.\ Sci.\ Univ.\ Tokyo Sect.\ IA, Math.\ \textbf{38} (1991), 251--298.

\bibitem[At]{atiyah} M. Atiyah,
\emph{Convexity and commuting hamiltonians}, Bull.\ London Math.\ Soc.\
\textbf{14} (1982), 1--15.

\bibitem[Au1]{audin:paper} M. Audin, 
\emph{Hamiltoniens p\'{e}riodiques sur les vari\'{e}t\'{e}s symplectiques
compactes de dimension 4}, G\'{e}om\'{e}trie symplectique et
m\'{e}canique, Proceedings 1988, C.\ Albert ed., Springer Lecture Notes in
Math.\ \textbf{1416} (1990).

\bibitem[Au2]{audin:book} M. Audin, 
\emph{The topology of torus actions on symplectic manifolds},
Progress in Mathematics \textbf{93}, Birkhäuser Verlag, Basel, 1991. 

\bibitem[BB]{b-b} A. Bialynicki-Birula,
\emph{Remarks on the action of an algebraic torus on $k^n$. II},
Bull.\ Acad.\ Polon.\ Sci.\ S\'er.\ Sci.\ Math.\ Astronom.\ Phys.\
\textbf{15}, (1967), 123--125.

\bibitem[BT]{BT} R. Bott and L.\ W. Tu, \emph{Differential forms in
algebraic topology}, Springer--Verlag, 1982.

\bibitem[De1]{De} T. Delzant,
\emph{Hamiltoniens p\'{e}riodques et image convexe de l'application
moment}, Bull.\ Soc.\ Math.\ France \textbf{116} (1988), 315--339.

\bibitem[De2]{Del2} T. Delzant,
\emph{Classification des actions hamiltoniennes compl\`ement
int\'egrables de rang deux}, Ann.\ Global Anal.\ Geom.\ \textbf{8}
(1990), 87--112.

\bibitem[DH]{DH} J. J. Duistermaat and G. J. Heckman,
\emph{On the variation in the cohomology of the symplectic form
of the reduced phase space}, Invent.\ Math.\ \textbf{69},
259--268 (1982).

\bibitem[FK]{FK} K.-H.\ Fieseler and L.\ Kaup,
\emph{On the geometry of affine algebraic $\C^*$-surfaces},
in: \emph{Problems in the theory of surfaces and their classification 
(Cortona, 1988)}, 111--140, Sympos.\ Math., XXXII, Academic Press, London,
1991.

\bibitem[F]{fintushel}  R.\ Fintushel,
\emph{Circle actions on simply connected 4-manifolds},
Trans.\ Amer.\ Math.\ Soc.\ \textbf{230} (1977), 147--171,
and \emph{Classification of circle actions on 4-manifolds},
Trans.\ Amer.\ Math.\ Soc.\ \textbf{242} (1978), 377--390.

\bibitem[GSj]{GSj}
V. Guillemin and R. Sjamaar, lecture at Mass.\ Inst.\ Tech., 1995.

\bibitem[GS1]{GS:convexity} V. Guillemin and S. Sternberg,
\emph{Convexity properties of the moment mapping, I},
Invent.\ Math.\ \textbf{67} (1982), 491--513.

\bibitem[GS2]{GS:normal} V. Guillemin and S. Sternberg,
\emph{A normal form for the moment map},
Differential geometric methods in mathematical physics,
(S. Sternberg, Ed.) Reidel, Dordrecht, Holland, 1984.

\bibitem[HS]{HS}
A. Haefliger  and  E. Salem, \emph{Actions of tori on orbifolds},
Ann.\ Global Anal.\ Geom.\ \textbf{9} (1991), 37--59.

\bibitem[I]{iglesias}
P.\ Iglesias, \emph{Les $\SO(3)$-vari\'et\'es symplectiques
et leurs classification en dimension 4},
Bull.\ Soc.\ Math.\ France \textbf{119} (1991), 371--396.

\bibitem[K1]{karshon:appendix}
Y.\ Karshon, 
\emph{Appendix} to \emph{Symplectic packings and algebraic geometry}
by D.\ McDuff and L.\ Polterovich, Invent.\ Math.\ \textbf{115},
431--434 (1994).

\bibitem[K2]{karshon:periodic}
Y.\ Karshon, \emph{Periodic Hamiltonian flows on four dimensional manifolds}, 
Mem.\ Amer.\ Math.\ Soc.\ \textbf{672} (1999).

\bibitem[KT]{gromov}
S.\ Tolman and Y.\ Karshon, in preparation.

\bibitem[Kn]{knop} F. Knop, private communication.

\bibitem[Ko]{kosinski} A. Kosinski, \emph{Differential manifolds},
Pure and Applied Mathematics \textbf{138}, Academic Press, Inc.,
Boston, MA (1993).

\bibitem[Kl]{koszul}
J.\ L.\ Koszul,
\emph{Sur certains groupes de transformations de Lie},
Colloque international du Centre National de la Recherche Scientifique
\textbf{52} (1953), 137--141.

\bibitem[LMTW]{LMTW} E. Lerman, E. Meinrenken, S. Tolman, and C. Woodward, 
\emph{Nonabelian convexity by symplectic cuts},
Topology \textbf{37} (1998), 245--259.

\bibitem[LV]{LV} D. Luna and Th.\ Vust,
\emph{Plongements d'espaces homog\`enes},
Comment.\ Math.\ Helv.\ \textbf{58} (1983), 186--245.

\bibitem[M]{marle}
C.\ M. Marle, { Mod\`{e}le d'action hamiltonienne d'un groupe de
Lie sur une vari\'{e}t\'{e} symplectique}, \emph{Rendiconti del
Seminario Matematico} \textbf{43} (1985), 227--251, Universit\`{a}e
Politechnico, Torino.

\bibitem[Mo]{moser} J. Moser,
\emph{On the volume elements on a manifold},
Trans.\ Amer.\ Math.\ Soc.\ \textbf{120} (1965), 286--294.

\bibitem[KKMS]{toroidal}
G.\ Kempf, F.\ F.\ Knudsen, D.\ Mumford, B.\ Saint-Donat,
\emph{Toroidal embeddings. I.}
Lecture Notes in Math.\ \textbf{339}, Springer-Verlag, 
Berlin-New York, 1973.

\bibitem[OR]{OR} P. Orlik and F. Raymond, 
"Actions of the torus on 4-manifolds." I and II, 
Trans.\ Amer.\ Math.\ Soc.\ \textbf{152} (1970), 531--559,
and  Topology \textbf{13} (1974), 89--12.

\bibitem[OW]{or-wag}   P.\ Orlik and P.\ Wagreich, 
\emph{Algebraic surfaces with $k^*$-actions}, 
Acta Math.\ \textbf{138} (1977), 43--81.

\bibitem[R]{rynes} J. Rynes,
\emph{Nonsingular affine $k^*$-surfaces},
Trans.\ Amer.\ Math.\ Soc.\ \textbf{332} (1992), 889--921.

\bibitem[Sch1]{schwarz:Topology} G.\ W. Schwarz,
\emph{Smooth functions invariant under the action of a compact Lie group},
Topology \textbf{14} (1975), 63--68.

\bibitem[Sch2]{schwarz:IHES} G.\ W. Schwarz,
\emph{Lifting smooth homotopies of orbit spaces},
I.H.E.S.\ Publ.\ Math.\ \textbf{51} (1980), 37--135.

\bibitem[T1]{T1}  D. A. Timash\"ev,
\emph{$G$-manifolds of complexity $1$}, 
(Russian) Uspekhi Mat.\ Nauk \textbf{51} (1996), no.\ 3 (309), 213--214;
English translation: Russian Math.\ Surveys \textbf{51} (1996), 567--568.

\bibitem[T2]{T2}  D. A. Timash\"ev,
\emph{Classification of $G$-manifolds of complexity $1$}, 
(Russian) Izv.\ Ross.\ Akad.\ Nauk Ser.\ Mat.\ \textbf{61} (1997), 127--162;
English translation: Izv.\ Math.\ \textbf{61} (1997), 363--397.

\bibitem[T]{traynor}
L.\ Traynor \emph{Symplectic packing constructions}, J.\ Diff.\ Geom.\ 
\textbf{42} (1995), 411--429.

\bibitem[W1]{w:lectures}
A. Weinstein, \emph{Lectures on Symplectic Manifolds},
CBMS Reg.\ conf.\ ser.\ in Math.\ \textbf{29}, 1977.

\bibitem[W]{woodward} C. Woodward, 
\emph{The classification of transversal multiplicity-free group actions},
Ann.\ Global Anal.\ Geom.\ \textbf{14} (1996), 3--42. 


\end{thebibliography}
\end{document}